\documentclass[12pt,leqno]{article}

\usepackage{amsthm}%,showidx}
\usepackage[usenames]{color}

\usepackage{amsmath}
\usepackage{amscd}
\usepackage{amsmath}
\usepackage{amssymb}
\usepackage{latexsym}
\makeatletter

\@addtoreset{equation}{section}
\makeatother

\newtheorem{thm}{Theorem}[section]
\newtheorem{pr}[thm]{Proposition}
\newtheorem{df}[thm]{Definition}
\newtheorem{lm}[thm]{Lemma}
\newtheorem{cor}[thm]{Corollary}

\input xy \xyoption{all} \CompileMatrices

\begin{document}

\title{Frobenius-Witt differentials
and regularity}
\author{Takeshi Saito}

\maketitle

\maketitle
\begin{abstract}
T.\ Dupuy, E.\ Katz, J.\ Rabinoff, D.\ Zureick-Brown
introduced the module
of total $p$-differentials
for a ring over 
${\mathbf Z}/p^2{\mathbf Z}$.
We study the same construction
for a ring over
${\mathbf Z}_{(p)}$
and prove a regularity criterion.
For a local ring,
the tensor product with
the residue field is constructed
in a different way
by O.\ Gabber, L.~Ramero.

In another article \cite{lex},
we use the sheaf of FW-differentials 
to define the cotangent bundle
and the
micro-support
of an \'etale sheaf.
\end{abstract}

Let $p$ be a prime number
and 
$P
=\dfrac{(X+Y)^p-X^p-Y^p}p
\in {\mathbf Z}[X,Y]$
be the polynomial appearing in
the definition of addition of
Witt vectors.
For a ring $A$ and
an $A$-module $M$,
we say a mapping $w\colon A\to M$
is a Frobenius-Witt derivation
(Definition \ref{dfFW})
or an FW-derivation for short
if for any $a,b\in A$, we have
\begin{align*}
w(a+b)\, &=
w(a)+
w(b)
-P(a,b)
\cdot w(p),\\
w(ab)\, &=
b^p\cdot w(a)+
a^p\cdot w(b).
\end{align*}
For rings
over ${\mathbf Z}/p^2{\mathbf Z}$,
such mappings are studied in
\cite{DKRZ} 
and called $p$-total derivation.
As we show in Lemma \ref{lmdel}.3,
we have $p\cdot w(a)=0$
for $a\in A$ if $A$ is a ring
over ${\mathbf Z}_{(p)}$
and then we may replace
$a^p, b^p$ in (\ref{eqLb})
by $F(\bar a),F(\bar b)$
for the absolute Frobenius morphism
$F\colon A/pA=A_1\to A_1$.
The equalities may be considered
as linearized variants
of 
those in the definition of
$p$-derivation 
\cite{Bu}
or equivalently
$\delta$-ring \cite{BS}.

After preparing basic properties
of FW-derivations in Section \ref{sFW},
we introduce the module
$F\Omega^1_A$ of 
FW-differentials for a ring $A$
endowed
with a universal FW-derivation
$w\colon A\to F\Omega^1_A$
in Lemma \ref{lmOm}.
If $A$ is a ring
over ${\mathbf Z}_{(p)}$,
then $F\Omega^1_A$
is an $A/pA$-modules
and the canonical morphism
$F\Omega^1_A
\to
F\Omega^1_{A/p^2A}$
is an isomorphism
by Corollary \ref{corA}.1.
Consequently, the generalization of
the definition does not introduce
new objects.
If $A$ itself is a ring over
${\mathbf F}_p$,
then the $A$-module 
$F\Omega^1_A$
is canonically identified with
the scalar extension 
$F^*\Omega^1_A$
of $\Omega^1_A$
by the absolute Frobenius
$F\colon A\to A$
by Corollary \ref{corA}.2.

For a local ring $A$ 
with residue field $k=A/{\mathfrak m}$
of characteristic $p$,
we show in Proposition \ref{prdx}
that the $k$-vector space $F\Omega^1_A
\otimes_Ak$
fits in an exact sequence
$0\to F^*({\mathfrak m}_A/
{\mathfrak m}_A^2)
\to
F\Omega^1_A
\otimes_Ak\to
F^*\Omega^1_k\to 0$
where $F^*$ denotes the
scalar extension by
the absolute Frobenius
$F\colon k\to k$.
We deduce from this in 
Corollary \ref{corGR} that
$F\Omega^1_A\otimes_Ak$
is canonically identified
with the $k^{1/p}$-vector space
${\bf \Omega}_A$
defined by Gabber and Ramero
in {\rm \cite[9.6.12]{GR}}
using an extension of $A$
involving the ring
of Witt vectors $W_2(k)$.
They use this module
to correct an incomplete proof
of a regularity criterion stated in
\cite[Chapitre 0, Th\'eor\`eme 22.5.4]{EGA4}.
In the case where $A$ is a discrete
valuation ring,
we construct injections from the
duals of the graded quotients
of the Galois groups
of Galois extensions
of the fraction field of $A$
by the filtration
of ramification groups
to twists of
$F\Omega^1_A
\otimes_Ak$
in \cite{gr}.

The main result is 
the following regularity criterion.
%and the openness assumption 
%on the regular part.
Under a suitable
finiteness condition,
we prove in Theorem \ref{thmreg}
that a noetherian local ring
$A$ with residue field
of characteristic $p$
is regular
if and only if
the $A/pA$-module 
$F\Omega^1_A$
is free of the correct rank,
using Proposition \ref{prdx}.
%of rank $\dim X$ on
%$X_{{\mathbf F}_p}$.

The construction of $F\Omega^1$
is sheafified and we obtain 
a sheaf of FW-differentials
$F\Omega^1_X$ on a scheme $X$.
We will use 
the sheaf of FW-differentials 
in \cite{lex}
to define the cotangent bundle
and the
micro-support
of an \'etale sheaf 
in mixed characteristic.
In the final section,
we study the relation of
$F\Omega^1_X$ with
${\cal H}_1$ of cotangent complexes.

The author thanks Luc Illusie for comments on
earlier versions,
for discussion on cotangent bundle
and on notation and terminology.
The author thanks Ofer Gabber
for indicating another construction of 
the module and 
for the reference to \cite{GR}
and \cite{DKRZ}.
The author thanks Alexander Beilinson 
for suggesting similarity 
to \cite{Bu}
and \cite{BS}.
The author thanks Akhil Mathew
heartily for pointing
out an error in Lemma \ref{lmWitt} and
Corollary \ref{corA}.3
and also for suggesting
an argument proving that in Theorem \ref{thmreg} the regularity condition (2) 
implies the flatness of $F\Omega^1_A$ without any finiteness assumption.

The research is partially supported by Grant-in-Aid (B) 19H01780.

\section{Frobenius-Witt derivation}\label{sFW}

We introduce Frobenius-Witt derivations
and study basic properties.

\begin{df}\label{dfFW}
Let $p$ be a prime number.

{\rm 1.}
Define a polynomial
$P\in {\mathbf Z}[X,Y]$
by
\begin{equation}
P=
\sum_{i=1}^{p-1}
\dfrac{(p-1)!}{i!(p-i)!}\cdot
X^iY^{p-i}.
\label{eqP}
\end{equation}

{\rm 2.}
Let $A$ be a ring
and $M$ be an $A$-module.
We say that a mapping
$w\colon A\to M$
is a Frobenius-Witt derivation
or FW-derivation for short
if the following condition is
satisfied:
For any $a,b\in A$, we have
\begin{align}
w(a+b)\, &=
w(a)+
w(b)
-P(a,b)
\cdot w(p),
\label{eqadd}\\
w(ab)\, &=
b^p\cdot w(a)+
a^p\cdot w(b).
\label{eqLb}
\end{align}
\end{df}

For a ring $A$ over ${\mathbf Z}_{(p)}$,
Definition \ref{dfFW}.2
is essentially the same
as \cite[Definition 2.1.1]{DKRZ}
since the condition (3) loc.~cit.~is
automatically satisfied
by Lemmas \ref{lmdel}.3 and  \ref{lmWitt}.2
below.

%\begin{lm}\label{lmPab}
%For every $a\in {\mathbf F}_p^\times$,
%there exists $b\in {\mathbf F}_p^\times$
%such that $P(a,b)\neq 0$.
%\end{lm}
%
%\proof{
%Define a polynomial
%$Q\in {\mathbf Z}_{(p)}[X,Y]$
%by $P=XYQ$.
%Then, since the monic polynomial
%$Q(a,Y)\in {\mathbf F}_p[Y]$
%is of degree $p-2$,
%there exists $b\in {\mathbf F}_p^\times$
%such that $Q(a,b)\neq 0$.
%\qed
%
%}
%\medskip

\begin{lm}\label{lmdel}
Let $A$ be a ring
and $w\colon A\to M$
be an FW-derivation.

{\rm 1.}
We have $w(1)=0$.
Let $a\in A$
and $n\in {\mathbf Z}$.
Then, we have
\begin{equation}
w(na)=
n\cdot
w(a)+a^p\cdot w(n).
\label{eqna}
\end{equation}
If $n\geqq 0$, we have
\begin{equation}
w(a^n)=
na^{p(n-1)}\cdot
w(a).
\label{eqan}
\end{equation}

{\rm 2.}
For $n\in {\mathbf Z}$,
we have
\begin{equation}
w(n)=
\dfrac{n-n^p}p\cdot w(p),
\label{eqn}
\end{equation}
In particular,
we have $w(0)=0$.

{\rm 3.}
Assume that
$A$ is a ring over ${\mathbf Z}_{(p)}$.
Then, for any $a\in A$,
we have $p\cdot w(a)=0$.
\end{lm}

In the most part of this
article, $A$ will be 
a ring over ${\mathbf Z}_{(p)}$.
Under this assumption, FW-derivations
$w\colon A\to M$
take values
in the $p$-torsion part of $M$
by Lemma \ref{lmdel}.3.

\proof{
1.
By putting $a=b=1$ in
(\ref{eqLb}),
we obtain $w(1)=0$.

Set $w_a(n)
=
n\cdot
w(a)+a^p\cdot w(n)$.
Then, by (\ref{eqadd})
and $P(n,m)a^p=P(na,ma)$,
we have
$w_a(n+m)
=w_a(n)+w_a(m)-
P(na,ma)\cdot w(p)$.
Since $w_a(1)=w(a)$,
we obtain (\ref{eqna})
by the ascending and
the descending inductions on $n$
starting from $n=1$ by (\ref{eqadd}).

For $n=0$, we have $w(a^0)=w(1)=0$.
By (\ref{eqLb}) and
induction on $n$,
we have
$w(a^{n+1})
=
a^pw(a^n)
+
a^{pn}w(a)
=
a^p\cdot na^{p(n-1)}w(a)
+
a^{pn}w(a)
=
(n+1)
a^{pn}w(a)$
and (\ref{eqan}) follows.

2.
Set $w_1(n)
=
\dfrac{n-n^p}p\cdot w(p)$.
Then, by binomial expansion,
$w_1$ satisfies
(\ref{eqadd}).
Hence we obtain (\ref{eqn})
similarly as in the proof of (\ref{eqna}).
By setting $n=0$
in (\ref{eqn}),
we obtain $w(0)=0$.

3.
Comparing 
(\ref{eqna}) and (\ref{eqLb}),
we obtain
$(n-n^p)\cdot w(a)=0$.
Since the $p$-adic valuation
$v_p(p-p^p)$ is $1$,
we obtain $p\cdot w(a)=0$.
\qed

}

\begin{lm}\label{lmWitt}
Assume that $A$ is flat
over ${\mathbf Z}$
and that the Frobenius
$F\colon A/pA\to A/pA$ is an isomorphism.

{\rm 1.}
The mapping
$w\colon A\to A/pA$
given by
$w(a^p+pb)\equiv b^p
\bmod pA$ for $a,b\in A$
is well-defined and
is an FW-derivation.

In particular, for $A=
{\mathbf Z}_{(p)}$,
the mapping 
$w\colon
{\mathbf Z}_{(p)}\to {\mathbf F}_p$
defined by
$w(a)=\dfrac{a-a^p}p\mod p$
is an FW-derivation.

{\rm 2.}
Let $\varphi\colon A\to A$
be an endomorphism satisfying
$\varphi(a)\equiv a^p\bmod p$
and let $\varphi_1\colon A\to A$
be the unique mapping
satisfying
$\varphi(a)=a^p+p\varphi_1(a)$.
Let 
$M$ be any $A$-module
and $w\colon A\to M$
be any FW-derivation.
Then, we have
$$w (r)
=\varphi_1(r)\cdot w(p)$$
for $r\in A$.
\end{lm}

\proof{
1.
Since 
$F\colon A/pA\to A/pA$ is 
assumed a surjection,
any element $r\in A$
may be written as
$r=a^p+pb$
for $a,b\in A$.
Since $(a+pb)^p\equiv a^p
\bmod p^2$,
the mapping $w$ is well-defined.
Since 
$$a^p+pb+
a'^p+pb'=(a+a')^p
+p(b+b'-P(a,a')),$$
we have $$w(a^p+pb+
a'^p+pb)
=
(b+b'-P(a,a'))^p
\equiv
w(a^p+pb)+
w(a'^p+pb)
-P(a^p+pb,
a'^p+pb')
\bmod p$$
and (\ref{eqadd}) is satisfied.
Since $$(a^p+pb)
(a'^p+pb')\equiv (aa')^p
+p(a'^pb+a^pb')
\bmod p^2,$$
we have $$w((a^p+pb)
(a'^p+pb))
=
(a'^pb+a^pb')^p\equiv
(a'^p+pb')^pw(a^p+pb)
+
(a^p+pb)^pw(a'^p+pb')
\bmod p$$
and (\ref{eqLb}) is satisfied.

For $a\in A={\mathbf Z}_{(p)}$,
we have $a=a^p+pb$ for $b\in {\mathbf Z}_{(p)}$
and $w(a)\equiv b^p\equiv b\bmod p$.
Alternatively, we can also verify directly
that the mapping
$w\colon {\mathbf Z}_{(p)}\to
{\mathbf F}_p$ defined by
$w(a)\equiv (a-a^p)/p
\bmod p$
satisfies (\ref{eqadd}) and (\ref{eqLb}).

2.
Since 
$F\colon A/pA\to A/pA$ is 
assumed a surjection,
we may write $r=a^p+pb$
for $a,b\in A$.
Since $\varphi (a)\equiv
r\bmod p$
implies
$\varphi (a)^p\equiv
r^p\bmod p^2$, we have
$\varphi(r)
=\varphi (a)^p+p
\varphi(b)
\equiv
r^p+pb^p
\bmod p^2$.
Further by (\ref{eqadd}),
 (\ref{eqan}),
(\ref{eqLb})
and by
$p\cdot w(p)=
p\cdot w(a)=
p\cdot w(b)=0$
in Lemma \ref{lmdel}.3,
we have
$w (r)=w(a^p)+w(pb)
=b^p\cdot w(p)
=\varphi_1(r)\cdot w(p)$.
\qed

}
\medskip

We give a relation between
FW-derivations
and 
Frobenius semi-linear derivations
for rings over ${\mathbf F}_p$.

\begin{lm}\label{lmBp}
Let $A$ be a ring,
$B$ be a ring over ${\mathbf F}_p$
and $g\colon A\to B$ be a morphism of
rings.
For a $B$-module $M$
and a mapping 
$w \colon A\to M$,
the following conditions
are equivalent:

{\rm (1)}
If we regard $M$ as
an $A$-module by
$g\colon A\to B$,
then
$w$ is
an FW-derivation
and $w(p)=0$.

{\rm (2)}
If we regard $M$ as
an $A$-module by the composition
$f=F\circ g\colon A\to B$
with the absolute Frobenius,
then
$w$ is a derivation.
\end{lm}

\proof{
(1)$\Rightarrow$(2):
If $w$ is an FW-derivation
satisfying $w(p)=0$,
then $w$ is additive by (\ref{eqadd}).
Further (\ref{eqLb})
means the Leibniz rule
with respect to the composition
$f=F\circ g\colon A\to B$.

(2)$\Rightarrow$(1):
If $w$ satisfies the Leibniz rule,
then we have $w(1)=1$.
Hence the additivity implies
$w(p)=0$ and (\ref{eqadd}).
The Leibniz rule
with respect to the composition
$f=F\circ g$
means (\ref{eqLb}) conversely.
\qed

}

\begin{lm}\label{lmAMI}
Let $A$ be a ring, $I\subset A$
be an ideal and
let $M$ be an $A$-module.
Then an FW-derivation
$w\colon A\to M$
induces an FW-derivation
$\bar w\colon
A/I\to M/(IM+A\cdot w(I))$.
\end{lm}

\proof{
By (\ref{eqadd}),
we have $w (a+b)
\equiv w(a)+w(b)
\bmod IM$ for $a\in A$ and
$b\in I$.
Hence $w$ induces a mapping
$\bar w\colon A/I\to M/(IM+A\cdot w(I))$.
Since $w$ satisfies 
(\ref{eqadd}) and (\ref{eqLb}),
$\bar w$ also satisfies 
(\ref{eqadd}) and (\ref{eqLb}).
\qed

}
\medskip

An extension of FW-derivation
to the ring of polynomials
is uniquely determined
by choosing the value
at the indeterminate.

\begin{pr}\label{prAX}
Let $A$ be a ring 
and $M$ be an $A[X]$-module.
Let $w\colon A\to M$
be an FW-derivation.

{\rm 1.}
Let $x\in M$ be
an element satisfying $px=0$.
Then, there exists
a unique FW-derivation
$\widetilde w\colon A[X]\to M$
extending $w$
and satisfying
$\widetilde w(X)=x$.

{\rm 2.}
If $A$ is a ring over ${\mathbf Z}_{(p)}$,
the mapping 
\begin{equation}
\{\text{FW-derivations
$\widetilde w\colon A[X]\to M$
extending $w$}\}
\to
M[p]=\{x\in M\mid p x=0\}
\label{eqMp}
\end{equation}
sending
$\widetilde w$ to $\widetilde w(X)$
is a bijection to
the $p$-torsion part of $M$.
\end{pr}

\proof{1.
For a polynomial
$f=\sum_{i=0}^na_iX^i\in A[X]$,
let
$f'\in A[X]$ denote the derivative
and set
\begin{align}
Q(f)\,&=
\sum_{\substack{
0\leqq
k_0,\ldots,k_n
<p,\\
k_0+\cdots+k_n=p}}
\dfrac{(p-1)!}
{k_0!\cdot k_1!\cdots k_n!}
\cdot
a_0^{k_0}(a_1X)^{k_1}
\cdots (a_nX^n)^{k_n}
\in A[X],
\label{eqQ0}\\
w^{(p)}(f)\,&=
\sum_{i=0}^nX^{pi}\cdot
w(a_i)
\in M.
\label{eqwp}
\end{align}
In {\rm (\ref{eqQ0})},
the summation is taken over the
integers 
$0\leqq k_0,\ldots,k_n<p$
satisfying
$k_0+\cdots+k_n=p$.

If $\widetilde w\colon
A[X]\to M$ is
an FW-derivation
extending $w$
and satisfying
$\widetilde w(X)=x$,
then 
by (\ref{eqadd}) and (\ref{eqLb})
we have 
\begin{equation}
\widetilde w(f)=
f'^p\cdot x
+w^{(p)}(f)-Q(f)\cdot w(p)
\label{eqdelx}
\end{equation}
for $f\in A[X]$.
Hence it suffices to
show that $\widetilde w$ 
defined by (\ref{eqdelx}) is actually
an FW-derivation.

For $f=\sum_{i=0}^na_iX^i,
\, g=\sum_{i=0}^nb_iX^i\in A[X]$,
set $$
f^{(p)}
=\sum_{i=0}^na_i^pX^{pi},
\qquad
R(f,g)
=\sum_{i=0}^nP(a_i,b_i)X^{pi}.$$
Then, we have
\begin{equation}
(f+g)^{(p)}=f^{(p)}+g^{(p)}+pR(f,g),\qquad
f^p=f^{(p)}+
pQ(f).
\label{eqfp}
\end{equation}
From this and
$(f+g)^p=f^p+g^p+pP(f,g)$,
by reducing to the universal case where
$A$ is flat over ${\mathbf Z}$,
we deduce
\begin{equation}
Q(f+g)=Q(f)+Q(g)+P(f,g)-R(f,g).
\label{eqQ}
\end{equation}

By (\ref{eqadd}),
we have
\begin{equation}
w^{(p)}(f+g)=
w^{(p)}(f)+w^{(p)}(g)-R(f,g)\cdot w(p).
\label{eqR}
\end{equation}
Since $px=0$, 
we have
$(f+g)'^p\cdot x=
f'^p\cdot x+
g'^p\cdot x$.
This and (\ref{eqR}) and
(\ref{eqQ}) show that
the mapping
$\widetilde w$ %\colon A[X]\to M$
satisfies (\ref{eqadd}).

We show that
the mapping
$\widetilde w$ %\colon A[X]\to M$
satisfies (\ref{eqLb}).
Since $px=0$,
we have
$(fg)'^p x=
f^p\cdot g'^px+
g^p\cdot f'^p x$.
Hence, we may assume $x=0$.
If $f$ and $g$ are monomials,
we have $Q(f)=Q(g)=Q(fg)=0$
and $w^{(p)}(fg)=f^p\cdot w^{(p)}(g)+
g^p\cdot w^{(p)}(f)$
and (\ref{eqLb}) is satisfied in this case.
For $f_1,f_2,g\in A[X]$,
we have
$w_0(f_1g+f_2g)
-
(w_0(f_1g)
+w_0(f_2g))
=
P(f_1g,f_2g)\cdot w(p)$
and
$(
(f_1+f_2)^p
w_0(g)
+g^pw_0(f_1+f_2))
-
(f_1^pw_0(g)
+g^pw_0(f_1)
+
f_2^pw_0(g)
+g^pw_0(f_2))
=
g^pP(f_1,f_2)\cdot w(p)$
by (\ref{eqR}) and (\ref{eqQ}).
Since $
P(f_1g,f_2g)=
g^pP(f_1,f_2)$,
the equality (\ref{eqLb}) 
follows by
induction on the numbers of
non-zero terms in $f$ and $g$.

2.
If $\widetilde w\colon
A[X]\to M$ is
an FW-derivation
extending $w$,
we have $\widetilde w(X)
\in M[p]$ by
the assumption that
$A$ is a ring over
${\mathbf Z}_{(p)}$
and Lemma \ref{lmdel}.3.
Thus, the assertion follows
from 1.
\qed

}

\section{Frobenius-Witt differentials}

We introduce the module
of Frobenius-Witt differentials
as the target of the universal
FW-derivation
and study basic properties.

\begin{lm}\label{lmOm}
Let $p$ be a prime number and
$A$ be a ring.
Then, there exists a universal pair
of an $A$-module
$F\Omega^1_A$
and an FW-derivation
$w\colon A
\to F\Omega^1_A$.
\end{lm}

\proof{
Let $A^{(A)}$ be
the free $A$-module
representing
the functor sending an
$A$-module $M$ to the set
${\rm Map}(A, M)$
and let $[\ ]\colon A
\to A^{(A)}$ denote
the universal mapping.
Define an $A$-module
$F\Omega^1_A$
to be the quotient of
$A^{(A)}$
by the submodule generated by
$[a+b]-[a]-[b]+P(a,b)[p]$
and $[ab]-a^p[b]-b^p[a]$
for $a,b\in A$.
Then, 
the pair of
$F\Omega^1_A$
and the composition
$w\colon A\to 
F\Omega^1_A$
of 
$[\ ]\colon A\to A^{(A)}$
with the canonical surjection
$A^{(A)}\to 
F\Omega^1_A$
satisfies the required universal property.
\qed

}
\medskip

\begin{df}
Let $A$ be a ring and 
$p$ be a prime number.
We call the $A$-module $F\Omega^1_A$
and $w\colon A\to 
F\Omega^1_A$
in Lemma {\rm \ref{lmOm}}
{\em the module of {\rm FW}-differentials of} $A$ and
{\em the universal {\rm FW}-derivation}.
For $a\in A$,
we call
$w(a)\in F\Omega^1_A$
{\em the {\rm FW}-differential} of $a$.
\end{df}
\medskip

If $A$ is a ring over ${\mathbf Z}_{(p)}$,
by Lemma \ref{lmdel}.3,
we have $p\cdot F\Omega^1_A=0$.
For a morphism $A\to B$ of rings,
the composition
$A\to B\to
F\Omega^1_B$ 
defines a canonical morphism
$F\Omega^1_A
\to
F\Omega^1_B$ 
and hence a $B$-linear morphism
\begin{equation}
F\Omega^1_A\otimes_AB
\to
F\Omega^1_B.
\label{eqFAB}
\end{equation}
%Define $F\Omega^1_{B/A}$
%to be its cokernel.

We study the module
of FW-differentials
of a quotient ring.

\begin{pr}\label{prAI}
Let $p$ be a prime number and
let $A$ be a ring.
Let $I\subset A$ be an ideal
and $B=A/I$ be the quotient ring.

{\rm 1.}
The canonical morphism
$F\Omega^1_A\otimes_AB
\to
F\Omega^1_B$
{\rm (\ref{eqFAB})}
induces an isomorphism
\begin{equation}
(F\Omega^1_A
\otimes_{A}B)/
(B\cdot w(I))
\to
F\Omega^1_B.
\label{eqAI}
\end{equation}
In particular, if the ideal $I$ is
generated by $a_1,\ldots, a_n\in A$,
we have an isomorphism
\begin{equation}
F\Omega^1_A
/(I\cdot F\Omega^1_A
+\sum_{i=1}^nA\cdot w(a_i))\to
F\Omega^1_B.
\label{eqai}
\end{equation}

{\rm 2.}
Let $B\to B'$ be
a morphism of rings
to a ring $B'$ over ${\mathbf F}_p$.
and let $F^*(I/I^2\otimes_BB')$
denote the tensor product
with respect to the absolute
Frobenius $F\colon B'\to B'$.
Then
the isomorphism {\rm (\ref{eqAI})}
defines an exact sequence
\begin{equation}
F^*(I/I^2\otimes_BB')
\to 
F\Omega^1_A\otimes_AB'
\to
F\Omega^1_B\otimes_BB'\to 0
\label{eqFI}
\end{equation}
of $B'$-modules.
\end{pr}

\proof{
1.
By Lemma \ref{lmAMI},
the universal FW-derivation
$w\colon A\to F\Omega^1_A$
induces an FW-derivation
$\bar w\colon B\to
M=
(F\Omega^1_A
\otimes_{A}B)/
(B\cdot w(I))$.
This defines a $B$-linear mapping
$F\Omega^1_B
\to M$
in the opposite direction.
Since the composition
$F\Omega^1_A\to F\Omega^1_B
\to M$ with the morphism
induced by $A\to B$
is the canonical surjection,
the composition
$M\to F\Omega^1_B
\to M$ with 
(\ref{eqAI})
is the identity of $M$.
Since the other composition 
$F\Omega^1_B\to
M\to F\Omega^1_B$
is also the identity,
(\ref{eqAI})
is an isomorphism.

If $I$ is
generated by $a_1,\ldots, a_n\in A$,
the image of
$w
\colon
I\otimes_{\mathbf Z}B
\to
F\Omega^1_A
\otimes_{A}B$
is generated by 
$w(a_1),\ldots,
w(a_n)$
as a $B$-module
by (\ref{eqadd})
and (\ref{eqLb}).

2.
The additive mapping
$w\colon I
\to 
F\Omega^1_A
\otimes_{A}B'$
is compatible with the
composition $A\to B'$ with 
the Frobenius
$F\colon B'\to B'$
by (\ref{eqLb}).
Hence $w$
induces a $B'$-linear mapping
$F^*(I/I^2\otimes_BB')
\to 
F\Omega^1_A
\otimes_{A}B'$.
Since its image is
$B'\cdot w(I)$,
the sequence (\ref{eqFI})
is exact by the isomorphism 
(\ref{eqAI}).
\qed

}

\begin{cor}\label{corA}
Let $A$ be a ring over
${\mathbf Z}_{(p)}$
and set $B=A/pA$
and $B_2=A/p^2A$.
For a $B$-module $M$,
let $F^*M$ denote the tensor product
$M\otimes_BB$ with respect
to the absolute Frobenius
$F\colon B\to B$.

{\rm 1.}
The $A$-module
$F\Omega^1_A$
is a $B$-module.
The morphism
$F\Omega^1_A
\to 
F\Omega^1_{B_2}$
induced by the surjection
$A\to B_2=A/p^2A$
is an isomorphism.

{\rm 2.}
The derivation
$d\colon A\to F^*\Omega^1_B$
is an FW-derivation
and defines an isomorphism
\begin{equation}
F\Omega^1_A
/(A\cdot\, w(p))
\to
F^*\Omega^1_B
\label{eqAB}
\end{equation}
of $B$-modules.
In particular, if
$p=0$ in $A=B$,
the isomorphism
{\rm (\ref{eqAB})}
gives an isomorphism
\begin{equation}
F\Omega^1_B
\to
F^*\Omega^1_B.
\label{eqB}
\end{equation}

{\rm 3.}
Assume that
$A$ is faithfully flat over ${\mathbf Z}_{(p)}$
and that
the Frobenius
$F\colon A/pA\to A/pA$ 
is an isomorphism.
Then,
$F\Omega^1_A$
is a non-zero $A/pA$-module
generated by $w(p)$.

In particular, if $A$ is a discrete valuation
ring with perfect residue field $k$
such that $p$ is a uniformizer,
then 
$F\Omega^1_A$
is a $k$-vector space of dimension $1$
generated by $w(p)$.

{\rm 4.}
Assume that $A$ is noetherian
and that the quotient
$A/\sqrt{pA}$
by the radical of the principal ideal
$pA$ is of finite type over
a field $k$ with finite $p$-basis.
%or is isomorphic to
%a localization of such a ring.
Then, the $A$-module
$F\Omega^1_A$
is of finite type.
\end{cor}

By Lemma \ref{lmdel}.3 and
Corollary \ref{corA}.1,
if $A$ is a ring over ${\mathbf Z}_{(p)}$,
an FW-derivation $w \colon A\to M$
is always induced by 
an FW-derivation $A/p^2A\to M[p]$
to the $p$-torsion part.
Examples after the proof
show that we cannot
relax the assumption in 4.~in
essential ways.

\proof{
1.
The $A$-module
$F\Omega^1_A$
is a $B$-module
by Lemma \ref{lmdel}.3.
Since $p\cdot F\Omega^1_A=0$,
we have
$w(p^2)=2p^p\cdot w(p)=0$.
Hence the isomorphism
$F\Omega^1_A/
(p^2\cdot F\Omega^1_A
+B_2\cdot\, w(p^2))
\to
F\Omega^1_{B_2}$
(\ref{eqai}) for $I=p^2A$
gives the required isomorphism
$F\Omega^1_A
\to 
F\Omega^1_{B_2}$.

2.
Let $M$ be a $B$-module.
By the universality of
$F\Omega^1_A$,
$A$-linear morphisms
$F\Omega^1_A/(
A\cdot w(p))\to M$
correspond bijectively to
FW-derivations
$w\colon A\to M$
satisfying $w(p)=0$.
By the universality of
$F^*\Omega^1_B$,
$B$-linear morphisms
$F^*\Omega^1_B\to M$
correspond bijectively to
usual derivations
$B\to M$
with respect to the Frobenius
$B\to B$.
Since $B=A/pA$,
usual derivations
$B\to M$ further
correspond bijectively to
derivations $A\to M$
with respect to the composition
$A\to B$ with the Frobenius.
Hence the assertion follows from
Lemma \ref{lmBp}.

3.
Since 
$F\colon A/pA\to A/pA$ is 
assumed a surjection,
we have $\Omega^1_{A/pA}=0$.
Hence by the isomorphism (\ref{eqAB}),
$F\Omega^1_A$
is an $A/pA$-module
generated by
one element $w(p)$.
Let
$w\colon A\to A/pA$
the FW-derivation
in Lemma \ref{lmWitt}.1
defined by
$w(a^p+pb)\equiv b^p
\bmod pA$ for $a,b\in A$.
If $A/pA\neq 0$, then we have
 $w(p)=1\neq 0$ and 
$F\Omega^1_A\neq 0$.

4.
A field $k$ is formally smooth
over ${\mathbf F}_p$
by \cite[Chapitre 0, Th\'eor\`eme (19.6.1)]{EGA4}.
Since the ideal
$\sqrt{pA}/pA
\subset A/pA=B$
is a nilpotent ideal
of finite type,
the morphism
$k\to A/\sqrt{pA}$
is lifted
to a morphism
$k\to A/pA=B$
of finite type.
Since $k$ is of finite $p$-basis,
the $k$-vector space
$\Omega^1_k$ is of finite dimension
and
the
$B$-module
$\Omega^1_B$
is of finite type
by the exact sequence
$\Omega^1_k\otimes_kB
\to \Omega^1_B
\to \Omega^1_{B/k}
\to 0$.
Thus, the assertion follows from
the isomorphism 
(\ref{eqAB})
of $B$-modules.
\qed

}

\medskip
\noindent{\it Example}
1.
Let $A=k$
be a field of characteristic $p>0$.
Then,
the $k$-vector space
$F\Omega^1_k
=F^*\Omega^1_k$
is finitely generated
if and only if $k$ 
has a finite $p$-basis.

2.
Let $k$ be a perfect field of characteristic
$p>0$ and let
$K\subset k((t))$
be a subextension
of finite type of transcendental
degree $n\geqq 1$ over $k$
as in \cite[Proposition 11.6]{GRa}.
Then, $A=k[[t]]\cap K\subset k((t))$
is a discrete valuation ring
with residue field $k$
and $\dim_kF\Omega^1_A
\otimes_Ak\leqq 1$
by (\ref{eqFI}).
Since the surjection
$A\to A/{\mathfrak m}_A^2
=k[t]/(t^2)$
induces a surjection
$F\Omega^1_A
\to F\Omega^1_{ A/{\mathfrak m}_A^2}
\neq 0$,
we have
$\dim_kF\Omega^1_A
\otimes_Ak=1$.
On the other hand,
we have $\dim_KF\Omega^1_A
\otimes_AK=
\dim_KF^*\Omega^1_K=n$.
Hence if $n>1$,
the $A$-module
$F\Omega^1_A$
is not finitely generated.
\medskip

Let $A\to B$ be a surjection
of rings over ${\mathbf Z}_{(p)}$
with kernel $I\subset A$.
Set $A_1=A/pA$ and $B_1=B/pB$
and let $I_1\subset A_1$ be the image of $I$.
Then the exact sequence (\ref{eqFI}),
the isomorphism (\ref{eqAB})
for $A$ and $B$
and 
the Frobenius pull-back
of the exact sequence
$I_1/I_1^2\to \Omega^1_{A_1}
\otimes_{A_1}B_1
\to \Omega^1_{B_1}\to 0$
define a commutative diagram
\begin{equation}
\begin{CD}
%@. B'\cdot w(p)@. B'\cdot w(p)@.\\
%@.@VV\cap V@VV\cap V@.\\
F^*(I/I^2\otimes_B B_1)
@>w>>F\Omega^1_A
\otimes_AB_1
@>>>F\Omega^1_B
\otimes_BB_1@>>>0
@.\\
@VVV@VVV@VVV@.\\
F^*(I_1/I_1^2)
@>>>F^*\Omega^1_{A_1}
\otimes_{A_1}B_1
@>>>F^*\Omega^1_{B_1}
@>>> 0
\\
\end{CD}
\end{equation}
of exact sequences.
The morphism
$w\colon 
F^*(I/I^2\otimes_B B_1)
\to F\Omega^1_A
\otimes_AB_1$
is induced by the
restriction of
the universal FW-derivation
$w\colon A\to F\Omega^1_A$
and the vertical arrows
are the canonical surjections.
By the isomorphism (\ref{eqAB}),
the bottom terms 
on the middle and right are
the quotients of the top terms
by the $B_1$-submodules
generated by $w(p)$.

\begin{pr}\label{prSA}
Let $p$ be a prime number and
let $A$ be a ring.

{\rm 1.}
If $A=\varinjlim_{\lambda\in \Lambda} A_\lambda$
is a filtered inductive limit,
the canonical morphism
$\varinjlim_{\lambda\in \Lambda}
F\Omega^1_{A_\lambda}
\to F\Omega^1_A$
is an isomorphism.

{\rm 2.}
Let $S\subset A$ be
a multiplicative subset.
Then, the canonical morphism
\begin{equation}
S^{-1}F\Omega^1_A
\to F\Omega^1_{S^{-1}A}
\label{eqS-1}
\end{equation}
is an isomorphism.

{\rm 3.}
Assume that
$A$ is a ring over
${\mathbf Z}_{(p)}$
and
let $B=A[X]$ be a polynomial ring.
Then, 
$F\Omega^1_B$
is the direct sum of 
$F\Omega^1_A
\otimes_AB$
with a free $B/pB$-module
of rank $1$ generated by
$w(X)$.
\end{pr}

\proof{
1.
For any $A$-module $M$, 
FW-derivations $A\to M$
are in bijection with
projective systems of
FW-derivations $A_\lambda\to M$.
$A$-linear mappings 
$\varinjlim_\lambda 
F\Omega^1_{A_\lambda}\to M$
are also in bijection with
projective systems of
$A_\lambda$-linear mappings 
$F\Omega^1_{A_\lambda}\to M$.
Hence the assertion follows
from the universality of
$F\Omega^1$.

2.
By (\ref{eqLb}),
the mapping 
$w\colon
S^{-1}A\to 
S^{-1}F\Omega^1_A$
given by
$w (a/s)=
1/s^p\cdot w(a)
-(a/s^2)^p\cdot w(s)$
is well-defined.
Since this is an FW-derivation,
we obtain a morphism
$F\Omega^1_{S^{-1}A}\to 
S^{-1}F\Omega^1_A$.
The composition
$F\Omega^1_A
\to F\Omega^1_{S^{-1}A}\to 
S^{-1}F\Omega^1_A$
is the canonical morphism
and the composition
$F\Omega^1_{S^{-1}A}\to 
S^{-1}F\Omega^1_A
\to F\Omega^1_{S^{-1}A}$
is the identity.
Hence the morphism
(\ref{eqS-1})
has an inverse and 
is an isomorphism.

3.
Let $M$ be a $B$-module.
Then, by Proposition \ref{prAX}
and by the universality of
$F\Omega^1$,
$B$-linear morphisms
$F\Omega^1_B
\to M$
corresponds bijectively
to pairs of
$A$-linear morphisms
$F\Omega^1_A
\to M$
and elements of $M[p]$.
Since these pairs
corresponds bijectively to
$B$-linear morphisms
$(F\Omega^1_A
\otimes_AB)
\oplus
(B/pB)
\to M$,
the assertion follows.
\qed

}

\medskip

We give a description
as an extension of
the fiber of the module
of FW-differentials
of a local ring
at the closed point.

\begin{pr}\label{prdx}
Let $A$ be a local ring
such that
the residue field $k=A/
{\mathfrak m}_A$
is of characteristic $p$.
For a $k$-vector space $M$,
let $F^*M$
denote the tensor product
$M\otimes_kk$
with respect to the Frobenius
$F\colon k\to k$.
Let
$w\colon 
F^*({\mathfrak m}_A/
{\mathfrak m}_A^2)
\to
F\Omega^1_A
\otimes_A k
=
F\Omega^1_A
/{\mathfrak m}_AF\Omega^1_A$
be the morphism
induced by the universal
FW-derivation
$w\colon A\to F\Omega^1_A$.
Then, 
the sequence 
\begin{equation}
\begin{CD}
0@>>>
F^*({\mathfrak m}_A/
{\mathfrak m}_A^2)
@>{w}>>
F\Omega^1_A
\otimes_A k
@>>>
F^*\Omega^1_k
@>>>0
\end{CD}
\label{eqdx}
\end{equation}
{\rm (\ref{eqFI})}
of $k$-vector spaces
is exact.
\end{pr}

\proof{
The exactness except the injectivity
of $w$ follows from (\ref{eqFI}).
First, we show the case where
$A$ is the localization
at a prime ideal
of a polynomial ring
$A_0=W_2(k)[T_1,\ldots,T_n]$ 
over the ring $W_2(k_0)$ of
Witt vectors of length 2
for a perfect field $k_0$
and an integer $n$.
Then, by Proposition \ref{prSA}.3
and \ref{prSA}.1
and Corollary \ref{corA}.1 and \ref{corA}.3,
the $A_0$-module
$F\Omega^1_{A_0}$
is free of rank $n+1$.
Hence 
by Proposition \ref{prSA}.2,
the $k$-vector space
$F\Omega^1_A
\otimes_A k
=F\Omega^1_{A_0}
\otimes_{A_0} k$
is of dimension $n+1$.

Let $d$ be the transcendence degree of
$k$ over $k_0$.
Then, we have
$\dim \Omega^1_k=d$.
The localization
$B$ at the inverse image
of ${\mathfrak m}_A$
by the composition
$W(k)[T_1,\ldots,T_n]\to
W_2(k)[T_1,\ldots,$
$T_n]
\to A$
is a regular local ring
of dimension $n+1-d$
and the canonical morphism
${\mathfrak m}_B/
{\mathfrak m}_B^2
\to
{\mathfrak m}_A/
{\mathfrak m}_A^2$
is an isomorphism.
Hence we have
$\dim{\mathfrak m}_A/
{\mathfrak m}_A^2=n+1-d$.
Since (\ref{eqdx})
is exact except possibly at
$F^*({\mathfrak m}_A/
{\mathfrak m}_A^2)$
by Proposition \ref{prAI}.2,
it follows that (\ref{eqdx}) is exact everywhere.

We show the general case.
By taking the limit,
we may assume that
$A$ is a localization of 
a ring $A_0$ of finite type over
${\mathbf Z}$.
By Corollary \ref{corA}.1,
we may assume that
$A_0$ is of finite type
over ${\mathbf Z}/p^2{\mathbf Z}
=W_2(k_0)$ for $k_0={\mathbf F}_p$.
We take a surjection
$B_0=W_2(k)[T_0,\ldots,T_n]
\to A_0$.
Let $B$ be the localization
of $B_0$ at the inverse image
of ${\mathfrak m}_A$
by the composition $B_0\to A_0\to A$
and let $I$ be the kernel of the
surjection $B\to A$.
Then, by Proposition \ref{prAI}.2,
we have a commutative diagram 
$$\begin{CD}
@.F^*(I\otimes_B k)
@=
F^*(I\otimes_B k)
@.@.\\
@.@VVV@VVV@.@.\\
0@>>>
F^*({\mathfrak m}_B/
{\mathfrak m}_B^2)
@>{w}>>
F\Omega^1_B
\otimes_B k
@>>>
F^*\Omega^1_k
@>>>0
\\
@.@VVV@VVV@|@.\\
@.
F^*({\mathfrak m}_A/
{\mathfrak m}_A^2)
@>{w}>>
F\Omega^1_A
\otimes_Ak
@>>>
F^*\Omega^1_k
@>>>0
\\
@.@VVV@VVV@.@.
\\
@.0@.0@.@.
\end{CD}$$
of exact sequences.
%Further by Proposition \ref{prAI}.3,
%the kernel of
%$F\Omega^1_B
%\otimes_B k
%\to
%F\Omega^1_A
%\otimes_Ak$
%equals the image
%of $F^*(I\otimes_B k)$
Hence the assertion follows.
\qed

}
\medskip

We prove a relation with 
${\bf \Omega}_A$ defined
by Gabber-Ramero.
For the definition of ${\bf \Omega}_A$,
we refer to \cite[9.6.12]{GR}.

\begin{cor}\label{corGR}
Let $A$ be a local ring 
such that the residue field
$k=A/{\mathfrak m}_A$
is of characteristic $p$.
Let ${\bf \Omega}_A$
be the $k^{1/p}$-vector space
defined in {\rm \cite[9.6.12]{GR}}
and 
regard $\mbox{\boldmath $d$}_A
\colon A\to {\bf \Omega}_A$
as an FW-derivation
by identifying the inclusion
$k\to k^{1/p}$
with the Frobenius $F\colon k\to k$.
Then, the morphism
$F\Omega^1_A\otimes_Ak
\to {\bf \Omega}_A$ 
induced by $\mbox{\boldmath $d$}_A$
is an isomorphism.
\end{cor}

\proof{
For a $k$-vector space $V$,
we identify $V\otimes_kk^{1/p}$
with $F^*V$ by identifying the inclusion
$k\to k^{1/p}$
with the Frobenius $F\colon k\to k$.
We consider the diagram
\begin{equation}
\begin{CD}
0@>>>
F^*({\mathfrak m}_A/
{\mathfrak m}_A^2)
@>>> F\Omega^1_A\otimes_Ak
@>>> F^*\Omega^1_{k/{\mathbf F}_p}
@>>>0
\\
@.@|@VVV@|@.\\
0@>>>
{\mathfrak m}_A/
{\mathfrak m}_A^2
\otimes_kk^{1/p}
@>>>
 {\bf \Omega}_A
@>>>
\Omega^1_{k/{\mathbf F}_p}
\otimes_kk^{1/p}
@>>>0.
\end{CD}
%\label{eqOmL}
\end{equation}
The upper line is exact
by Proposition \ref{prdx}
and the lower exact sequence
is defined in
\cite[Proposition 9.6.14]{GR}.
The middle vertical arrow is
induced by the FW-derivation
$\mbox{\boldmath $d$}_A
\colon A\to {\bf \Omega}_A$
and the diagram is commutative.
Hence the assertion follows.
\qed

}
\medskip

We give a criterion of regularity
which will be used
in the proof of the main theorem
in the next section.

\begin{cor}\label{corXZx}
Let $A$ be a regular
local ring
such that the residue field
$k=A/{\mathfrak m}_A$
is of characteristic $p$.
Let $B=A/I$ be 
the quotient by an ideal 
$I\subset {\mathfrak m}_A$.
We set $A_1=A/pA$,
$B_1=B/pB$,
and 
for a $B_1$-module $M$,
let $F^*M=M\otimes_{B_1}
B_1$ denote
the tensor product with
respect to the Frobenius $F\colon 
B_1\to B_1$.
Let $w\colon
F^*(I\otimes_AB_1)
\to
F\Omega^1_A\otimes_A
B_1$
be the morphism
induced by the universal
FW-derivation
$w\colon A\to F\Omega^1_A$.

We consider the following conditions:

{\rm (1)} The sequence 
\begin{equation}
0\to
F^*(I\otimes_AB_1)
\overset w\longrightarrow
F\Omega^1_A\otimes_A
B_1\longrightarrow
F\Omega^1_B
\to
0
\label{eqXZx}
\end{equation}
of $B_1$-modules
is a split exact sequence.

{\rm (2)}
$B$ is regular.

{\rm 1.}
We always have
{\rm (1)}$\Rightarrow${\rm (2)}.

{\rm 2.}
Assume that
$F\Omega^1_A$
is a free $A_1$-module of
finite rank.
Then, we have
{\rm (2)}$\Rightarrow${\rm (1)}
and
$F\Omega^1_B$
is a free $B_1$-module of
finite rank.
\end{cor}

\proof{
First, we show that the condition
(2) is equivalent to the following condition:

(2$'$)
The sequence
$0\to I\otimes_A k
\to {\mathfrak m}_A
/{\mathfrak m}_A^2
\to
{\mathfrak m}_B/
{\mathfrak m}_B^2
\to 0$
is exact.

(2)$\Rightarrow$(2$'$):
The condition (2) means that
$I$ is generated by a part of
regular system of parameters
of $A$
by \cite[Chapitre 0, Corollaire (17.1.9)]{EGA4}.
This condition
means that
the images of
a minimal system of generators
of $I$ form
a basis of the kernel of 
$ {\mathfrak m}_A
/{\mathfrak m}_A^2
\to
{\mathfrak m}_B/
{\mathfrak m}_B^2$.
Hence the condition (2)
implies (2$'$).

(2$'$)$\Rightarrow$(2):
Conversely,
a lifting of the basis
of $I\otimes_A k$
is a part of
regular system of parameters
of $A$
and is a system of generators of $I$
by Nakayama's lemma.

By Proposition \ref{prdx} for $A$ and $B$,
(2$'$) is equivalent to the following:

(1$'$)
The sequence
\begin{equation}
0\to
F^*( I\otimes_A k)
\overset w\longrightarrow
F\Omega^1_A\otimes_Ak
\longrightarrow
F\Omega^1_B\otimes_Bk
\to
0
\end{equation}
induced by (\ref{eqXZx}) is exact.

1.
The condition (1) obviously implies (1$'$).

2. 
Since $F^*(I\otimes_AB_1)$
and
$F\Omega^1_A\otimes_AB_1$
are free $B_1$-modules
of finite rank,
the condition (1$'$) conversely implies (1)
and that
$F\Omega^1_B$
is a free $B_1$-module
of finite rank.
\qed

}

\begin{lm}\label{lmAB}
Let $f\colon A\to B$
be a morphism\
%of relative dimension $n$
of rings over ${\mathbf Z}_{(p)}$
and set
$A_1=A/pA$
and $B_1=B/pB$.
Then, the isomorphism {\rm (\ref{eqAB})}
induces an isomorphism
\begin{equation}
{\rm Coker}
(F\Omega^1_A
\otimes_AB\to F\Omega^1_B)
\to
F^*\Omega^1_{B_1/A_1}.
\label{eqXY}
\end{equation}
\end{lm}

\proof{
By the isomorphism (\ref{eqAB})
for $A$ and $B$
and its functoriality,
we have a commutative
diagram 
$$\begin{CD}
B_1
@>{\cdot w(p)}>>
F\Omega^1_A
\otimes_{A_1}B_1
@>>>
F^*(\Omega^1_{A_1}
\otimes_{A_1}B_1)
@>>>
0\\
@|@VVV@VVV@.\\
B_1@>{\cdot w(p)}>>
F\Omega^1_B
@>>>
F^*\Omega^1_{B_1}
@>>>0
\end{CD}$$
of exact sequences and the assertion follows.
\qed

}
\medskip

We give a criterion
for the smoothness.

\begin{pr}\label{prsm}
Let $f\colon A\to B$
be a morphism of finite presentation
%of relative dimension $n$
of rings over ${\mathbf Z}_{(p)}$
and set
$A_1=A/pA$
and $B_1=B/pB$.
We consider the sequence
\begin{equation}
\begin{CD}
0@>>>
F\Omega^1_A\otimes_AB
@>{\rm (\ref{eqFAB})}>>
F\Omega^1_B
@>>>
F^*\Omega^1_{B_1/A_1}
@>>>
0
\end{CD}
\label{eqXYFx}
\end{equation}
of $B_1$-modules

{\rm 1.}
Assume that $f$ is smooth.
Then, the sequence
{\rm (\ref{eqXYFx})}
is a split exact sequence and
{\rm (\ref{eqXY})} is an isomorphism of 
projective $B_1$-modules of finite rank.

{\rm 2.}
Let ${\mathfrak q}$ be
a prime ideal of $B$
such that 
the residue field $k=B_{\mathfrak q}
/{\mathfrak q}B_{\mathfrak q}$
is of characteristic $p$
and let ${\mathfrak p}\subset A$
be the inverse image of ${\mathfrak q}$.
Assume that
$A_{\mathfrak p}$ and 
$B_{\mathfrak q}$ are regular and that 
{\rm (\ref{eqXYFx})}
is a split exact sequence
after $\otimes_BB_{\mathfrak q}$.
Then $f\colon A\to B$
is smooth at ${\mathfrak q}$.
\end{pr}

\proof{
1.
Since $f$ is smooth,
the $B_1$-module 
$F^*\Omega^1_{B_1/A_1}=
{\rm Coker}
(F\Omega^1_A
\otimes_AB\to F\Omega^1_B)$ is
projective of finite rank.

If $B=A[T]$,
the assertion follows from
Proposition \ref{prSA}.3.
Since the question
is local on ${\rm Spec}\, B$, it suffices
to show that 
the morphism (\ref{eqFAB})
is an isomorphism
assuming that $A\to B$ 
is \'etale.

Since $A\to B$ is \'etale,
after a localization,
there exists 
a monic polynomial
$f\in A[T]$
such that ${\rm Spec}\, B$
is isomorphic to an open 
subscheme of 
${\rm Spec}\, A[T]/(f)[1/f']$
by \cite[Th\'eor\`eme (18.4.6)]{EGA4}.
Hence
we may further assume 
$B=A[T]/(f)[1/f']$
for a monic polynomial
$f\in A[T]$.
Then, by Proposition \ref{prSA}.3
and \ref{prSA}.2
and Proposition \ref{prAI}.1,
the $B/pB$-module
$F\Omega^1_B$
is the quotient of
$(F\Omega^1_A
\otimes_AB)
\oplus (B/pB\cdot w(T))$
by the submodule
generated by 
$\widetilde w(f)=
f'^{(p)}(T^p)\cdot w(T)
+w^{(p)}(f) +
Q(f)\cdot w(p)$
in the notation of
the proof of Proposition \ref{prAX}.
Since $f'^{(p)}(T^p)
\equiv f'^p\bmod pB$
is invertible in $B/pB$
and $w^{(p)}(f) +
Q(f)\cdot w(p)
\in F\Omega^1_A
\otimes_AB$,
the morphism
$F\Omega^1_A
\otimes_AB
\to\bigl(
(F\Omega^1_A
\otimes_AB)
\oplus (B/pB\cdot w(T))
\bigr)
/B\cdot \widetilde w(f)$
is an isomorphism
as required.

2.
Since the assertion is local
by Proposition \ref{prSA}.2,
we may assume that 
$A=A_{\mathfrak p}$.
We take a surjection
$C=A[T_1,\ldots,T_n]\to B$
and let $C_{\mathfrak r}$ be the localization
at the inverse image ${\mathfrak r}$ of 
${\mathfrak q}$.
Then, we have a split exact sequence
\begin{equation}
0\to
F\Omega^1_A\otimes_AC
\to
F\Omega^1_{C}
\to
F^*(\Omega^1_{C/A}\otimes_{C}
C/pC)
\to 
0
\label{eqABC}
\end{equation}
by Proposition \ref{prSA}.3.

By Proposition \ref{prdx}
for $C_{\mathfrak r}$ and 
$B_{\mathfrak q}$,
we have a commutative diagram
$$\begin{CD}
0
@>>>
F^*
({\mathfrak r}C_{\mathfrak r}/
{\mathfrak r}^2C_{\mathfrak r})
@>>>
F\Omega^1_C\otimes_Ck
@>>>
F^*\Omega^1_k
@>>>0\\
@.@VVV@VVV@|@.\\
0@>>>
F^*({\mathfrak q}B_{\mathfrak q}/
{\mathfrak q}^2B_{\mathfrak q})
@>>>
F\Omega^1_B
\otimes_Bk
@>>>
F^*\Omega^1_k
@>>>0
\end{CD}
$$
of exact sequences.
The vertical arrows
are surjections.
Since the kernel $I$
of the surjection
$C_{\mathfrak r}\to 
B_{\mathfrak q}$ of regular local rings
is generated by a part
of a regular system of
local parameters,
the sequence
$0\to I
\otimes_{C_{\mathfrak r}}k
\to
{\mathfrak r}C_{\mathfrak r}/
{\mathfrak r}^2C_{\mathfrak r}
\to
{\mathfrak q}B_{\mathfrak q}/
{\mathfrak q}^2B_{\mathfrak q}
\to 0$
is exact.
Hence 
we obtain an exact sequence
\begin{equation}
0\to F^*(I
\otimes_{C_{\mathfrak r}}k)
\to 
F\Omega^1_C\otimes_Ck
\to 
F\Omega^1_B
\otimes_Bk
\to 0.
\label{eqBC}
\end{equation}

If 
$
F\Omega^1_A\otimes_AB_{\mathfrak q}
\to
F\Omega^1_{B_{\mathfrak q}}$
is a split injection,
by (\ref{eqABC})
and (\ref{eqBC}),
the induced morphism
$F^*(I\otimes_{C_{\mathfrak r}}k)
\to
F^*(\Omega^1_{C/A}\otimes_Ck)$
is an injection.
This means that
the morphism 
$I/I^2
\to 
\Omega^1_{C_{\mathfrak r}/A}
\otimes_{C_{\mathfrak r}}
B_{\mathfrak q}$
of free $B_{\mathfrak q}$-modules
is a split injection.
Since $A\to C$ is smooth,
$A\to B$ is also smooth
at ${\mathfrak q}$.
\qed

}

\section{Regularity criterion}

We recall some facts
from commutative algebra
and field theory
in positive characteristic
used in the proof of
the main theorem.
Let $W$ be 
a noetherian complete
local ring and assume
that the characteristic
of the residue field is
a prime number $p$.
Then, $W$ is
said to be a Cohen ring
\cite[Chapitre 0, D\'efinition (19.8.4)]{EGA4}
if $W$ is flat over
${\mathbf Z}_p$
and $W/pW$
is a field,
or equivalently if
$W$ is an absolutely
unramified discrete valuation ring.

\begin{thm}\label{thmCohen}
Let $p$ be a prime number.

{\rm 1. (\cite[Chapitre 0, Th\'eor\`eme (19.8.2) (i)]{EGA4})}
Let $W$ be a Cohen ring
such that the residue field
is a field of characteristic $p$.
Then $W$ 
is formally smooth over
${\mathbf Z}_p$.

{\rm 2. (\cite[Chapitre 0, Th\'eor\`eme (19.8.6) (ii)]{EGA4})}
If $k$ is a field of characteristic $p$,
there exists a Cohen ring
$W$ such that
the residue field is isomorphic to
$k$.
\end{thm}

A local noetherian ring $A$
is said to be of complete
intersection if
its completion
$\hat A$
is isomorphic to
the quotient
of a regular complete local
noetherian ring $B$
by the ideal generated
by a regular sequence of $B$
\cite[Chapitre IV, D\'efinition (19.3.1)]{EGA4}.
Let $f\colon X\to S$
be a flat morphism of finite
type of noetherian schemes
and $x\in X, s=f(x)\in S$.
We say that $X$
is locally of complete intersection
relatively to $S$ at $x$
if the local ring ${\cal O}_{X_s,x}$
of the fiber $X_s=X\times_Ss$
is of complete intersection
\cite[Chapitre IV, D\'efinition (19.3.6)]{EGA4}.
Let $i\colon X\to Y$ be a closed
immersion of schemes
of finite type over a noetherian
schemes $S$
and $x\in X$.
We say that $i$ is transversally
regular relatively to $S$
at $x$ if on a neighborhood
$V\subset Y$ of $x$
there exists a regular sequence
$(f_i; 1\leqq i\leqq n)$
generating
the ideal ${\cal I}_X\subset
{\cal O}_Y$ defining $X$
such that 
${\cal O}_Y/(f_i; 1\leqq i\leqq j)$
are flat over $S$
for $1\leqq j\leqq n$
\cite[Chapitre IV, D\'efinition (19.2.2)]{EGA4}.

\begin{pr}\label{prci}
{\rm 1.
(\cite[Chapitre IV, Proposition (19.3.2)]{EGA4})}
Let $A=B/I$ be a quotient ring
of a regular local noetherian ring $B$.
Then, $A$ is of complete intersection
if and only if $I$ is generated
by a regular sequence of $B$.

{\rm 2
(\cite[Chapitre IV, Proposition (19.3.7)]{EGA4})}
Let $i\colon X\to Y$
be a closed immersion
of flat schemes of finite
type over a noetherian scheme $S$
and $x\in X$.
Assume that $Y$ is smooth over $S$.
Then, the immersion 
 $i$ is transversally
regular relatively to $S$
at $x$
if and only if
$X$ is locally of complete intersection
relatively to $S$ at $x$.
\end{pr}

\begin{thm}\label{thmpba}
Let $k$ be a field of characteristic $p>0$.

{\rm 1.
(\cite[Section 13, No.~2, Th\'eor\`eme 2 c)]{B})}
If $[k:k^p]=n$ is finite,
$\dim_k\Omega^1_{k/{\mathbf F}_p}=n$.

{\rm 2.
(\cite[Section 16, No.~6, Corollaire 3]{B})}
Let $k_1$ be a subfield such that
$k$ is finitely generated over $k_1$
of transcendental degree $d$
and that $[k_1:k_1^p]$ is finite.
Then 
$[k:k^p]=p^d\cdot [k_1:k_1^p]$.
\end{thm}

We say that a local ring $A$
is essentially of finite type
over a field $k$
if $A$ is isomorphic 
to the localization at
a prime ideal of
a ring of finite type
over $k$.
We state and prove
the regularity criterion.

\begin{thm}\label{thmreg}
Let $A$ be a noetherian local ring
with residue field
$k=A/{\mathfrak m}_A$
of characteristic $p$.
Assume that $k$ has a finite $p$-basis
and set  $d=\dim A$,
$[k:k^p]=p^r$ 
and $A_1=A/pA$.
We consider the following conditions:

{\rm (1)}
The $A_1$-module
$F\Omega^1_A$
is free of rank $d+r$.

{\rm (1$'$)}
The $k$-vector space
$F\Omega^1_A
\otimes_Ak$
is of dimension $d+r$.

{\rm (2)}
%$X$ is regular %of dimension $n$
%on a neighborhood of
%$X_{{\mathbf F}_p}$.
%
%{\rm (1$'$)}
$A$ is regular. 
%of dimension $n-d$
%for $d=\dim\overline{\{x\}}$.

{\rm 1.}
We always have
{\rm (1)}$\Rightarrow${\rm (1$'$)}$\Rightarrow${\rm (2)}.

{\rm 2.}
Assume that the quotient
$A/\sqrt{pA}$ by the radical
of the principal ideal $pA$
is essentially 
of finite type over
a field $k_1$ with finite $p$-basis
and that either of
the following conditions 
is satisfied:

{\rm (a)}
$A$ is flat over ${\mathbf Z}_{(p)}$.

{\rm (b)}
$A$ is a ring over ${\mathbf F}_p$.

\noindent
Then the 3 conditions
are equivalent.
%
%\noindent
%Then, 
%we have implications
%{\rm (1)}$\Rightarrow${\rm (2)}$\Rightarrow${\rm (2$'$)}$\Rightarrow${\rm (1$'$)}.
%If the subset
%${\rm Reg}(X)
%\subset X$ consisting
%of regular points is an open subset,
%the conditions are equivalent
%to each other.
\end{thm}

Let $A$ be the discrete valuation ring
in Example 2 after Corollary \ref{corA}.
Then $A$ satisfies (2)
and (1$'$) for $d=1$, $r=0$
but not (1)
unless $n=1$.

\proof{
1.
The implication
(1)$\Rightarrow$(1$'$) is obvious.
We show 
(1$'$)$\Rightarrow$(2).
By Proposition \ref{prdx},
we have
$\dim_k {\mathfrak m}_A/{\mathfrak m}^2_A
=
\dim_k F\Omega^1_A
\otimes_Ak
-
\dim_k\Omega^1_k
=
(d+r)-r=d
=\dim A$.
Hence $A$ is regular.

2.
It suffices to show (2)$\Rightarrow$(1).
First, we show the case (a).
Assume that $A$ is flat over ${\mathbf Z}_{(p)}$.
Let $W$ be a Cohen ring
with residue field $k_1$.
Then, since $W_2=W/p^2W$
is formally smooth over
${\mathbf Z}/p^2{\mathbf Z}$
by Theorem \ref{thmCohen}.2
and the ideal $\sqrt{pA}/p^2A
\subset  A_2=A/p^2A$
is nilpotent,
the morphism
$k_1\to A/\sqrt{pA}$
is lifted to a morphism
$W_2\to A_2$.
By the exact sequence
$0\to A/pA
\to A/p^2A
\to A/pA\to 0$,
we have
${\rm Tor}_1^{W_2}(A_2,k_1)=0$
and
the ring $A_2$ is flat over $W_2$.

Since the ideal $\sqrt{pA}/p^2A
\subset A_2$
is finitely generated,
there exists a morphism
$C_2=W_2[T_1,\ldots,T_N]
\to A_2$ over $W_2$
for an integer $N\geqq0$
such that 
for the localization $B_2$
of $C_2$ at the inverse image
of ${\mathfrak m}_{A_2}$,
the induced morphism
$B_2\to A/\sqrt{pA}$
is a surjection and
that the image
$C_2\to A_2$
contains a system of
generators of
$\sqrt{pA}/p^2A
\subset A_2$.
Then, since
$\sqrt{pA}/p^2A$ is nilpotent,
the local morphism
$B_2\to A_2$ is a surjection.

%The rest of the proof is similar to
%that of Corollary \ref{corXZx}.2. 
Set $B_1=B_2/pB_2$,
$C_1=C_2/pC_2$
and $n=d+{\rm tr.\ deg}_{k_1}k$.
Since $B_1$ is the local ring
of $k_1[T_1,\ldots,T_N]$
at a prime ideal with the residue field $k$,
we have $\dim B_1
=N-{\rm tr.\ deg}_{k_1}k$.
Since $A$ is regular
and $p\in A$
is a non-zero divisor,
the quotient $A_1=A/pA$ is of complete
intersection.
Since $B_1$ is regular,
the kernel $I_1$ of
the surjection
$B_1\to A_1$
is generated by a regular sequence 
of length $\dim B_1-(\dim A-1)
=(N-{\rm tr.\ deg}_{k_1}k)-(d-1)
=N-n+1$.

Let $X\subset {\mathbf A}^N_{W_2}
={\rm Spec}\, W_2[T_1,\ldots,T_N]$
be a closed subscheme
such that $A_2$ is isomorphic
to the local ring at a point $x\in X$.
Since $A_2$ is flat over $W_2$
and  $A_1$ is of complete
intersection,
the closed immersion
$X\to {\mathbf A}^N_{W_2}$
is transversally regular
relatively to $W_2$
at $x$
by Proposition \ref{prci}.2.
Hence  
the kernel $I_2$ of
the surjection
$B_2\to A_2$
is also generated by a regular sequence 
of length $N-n+1$
and the canonical surjection
$I_2/I_2\otimes_{A_2}A_1
\to 
I_1/I_1$
is an isomorphism
of free $A_1$-modules
of rank $N-n+1$.

The canonical
morphism $F\Omega^1_A
\to F\Omega^1_{A_2}$
is an isomorphism 
of $A_1$-modules by Corollary \ref{corA}.1.
Hence, we obtain an
exact sequence
\begin{equation}
F^*(I_2/I_2\otimes_{A_2}A_1)
\to F\Omega^1_{C_2}
\otimes_{C_1}A_1
\to F\Omega^1_A
\to 0
\label{eqres}
\end{equation}
of $A_1$-modules
by Proposition \ref{prAI}.2
and $F^*(I_2/I_2\otimes_{A_2}A_1)
=F^*(I_1/I_1^2)$
is a free $A_1$-module
of rank $N-n+1$.

%We compute the dimensions
%of the $k$-vector spaces
%obtained by taking tensor product
%of the $A_1$-modules
%in (\ref{eqres}).
Set
$[k_1:k_1^p]=p^{r_1}$.
We have
$\dim_{k_1}\Omega^1_{k_1}=r_1$
by Theorem \ref{thmpba}.1.
The $W_2$-module
$F\Omega^1_{W_2}$
is a $k_1$-vector space
by Corollary \ref{corA}.1
and is of dimension $r_1+1$
by Proposition \ref{prdx}.
Hence by Proposition \ref{prSA}.3,
the $C_2$-module
$F\Omega^1_{C_2}$
is a free $C_1$-module
of rank $N+r_1+1$.

%By  Proposition \ref{prdx},
We have
$r=\dim_k\Omega^1_k=
\dim_{k_1}\Omega^1_{k_1}
+{\rm tr.\, deg}_{k_1}k$
by Theorem \ref{thmpba}.
Since $A$ is regular,
by Proposition \ref{prdx},
the $k$-vector space
$F\Omega^1_A
\otimes_Ak$
is of dimension $d+r=
d+
{\rm tr.\, deg}_{k_1}k+r_1=n+r_1$.

Since $N+r_1+1=
(N-n+1)+(n+r_1)$,
the exact sequence
(\ref{eqres}) induces 
an exact sequence
$0\to F^*(I_1/I_1^2)
\otimes_{A_1}k
\to F\Omega^1_{C_2}
\otimes_{C_1}k
\to F\Omega^1_A
\otimes_{A_1}k
\to 0$.
Consequently
the morphism
$F^*(I_1/I_1^2)
\to F\Omega^1_{C_2}
\otimes_{C_1}A_1$
of free $A_1$-modules of finite rank
is a split injection
and $F\Omega^1_A$
is a free $A_1$-module
of rank $d+r$.

The proof in the case (b)
is similar and easier.
Since $k$ is formally smooth
over ${\mathbf F}_p$, we may
assume that
$A$ is the localization
at a prime ideal
of a ring $B$ of finite
type over $k_1$
and take a surjection
$C=k_1[T_1,\ldots,T_N]\to B$.
%Then, the immersion
%is regular of codimension $N-n$
%and we have an exact sequence
%\begin{equation}
%F^*N_{X/P}
%\to F\Omega^1_P
%\otimes_{{\cal O}_P}{\cal O}_X
%\to F\Omega^1_X
%\to 0
%\label{eqres2}
%\end{equation}
%of ${\cal O}_X$-modules.
By Corollary \ref{corA}.2,
$F\Omega^1_C$
is isomorphic
to the free $C$-module
$F^*\Omega^1_C$
of rank $N+r_1$.
Hence it suffices to apply 
Corollary \ref{corXZx}.2
to the localization of $C\to A$.
%of dimension $n-d$.
%
%If ${\rm Reg}(X)
%\subset X$ is an open subset,
%the conditions {\rm (1)}
%and {\rm (1$'$)} are equivalent.
\qed

}

\begin{cor}\label{corXZ}
Let $A\to A/I=B$ be a surjection
of regular local rings.
Assume that the quotient
$A/\sqrt{pA}$ by the radical
of the principal ideal $pA$
is essentially of finite type over
a field $k_1$ with finite $p$-basis.
Then for
$B_1=B/pB$,
the sequence 
\begin{equation}
0\to
F^*(I/(I^2+pI))
\overset w\longrightarrow
F\Omega^1_A\otimes_A B_1
\longrightarrow
F\Omega^1_B
\to
0
\label{eqXZ}
\end{equation}
of 
$B_1$-modules
is a split exact sequence.
%
%{\rm 2.}
%Let $x\in X
%_{{\mathbf F}_p}$.
%If 
%${\cal O}_{X,x}$ is regular
%and if the stalk of
%{\rm (\ref{eqXZ})} at $x$
%is a split exact sequence,
%the local ring
%${\cal O}_{Z,x}$ is regular.
\end{cor}

\proof{
%Let $[k:k^p]=p^r$.
%By Theorem \ref{thmreg} 
%(1)$\Rightarrow$(2),
Since the $A/pA$-module
$F\Omega^1_A$
%\otimes_{{\cal O}_X}
%{\cal O}_{Z_{{\mathbf F}_p}}$
%and
%$F\Omega^1_Z\otimes_{{\cal O}_Z}
%{\cal O}_{Z_{{\mathbf F}_p}}$
is free of finite rank
by Theorem \ref{thmreg}.2,
the assertion follows from Corollary \ref{corXZx}.2. 
%$\dim X+r$ and $\dim Z+r$ respectively.
%Since 
%$
%F^*(N_{Z/X}\otimes_{{\cal O}_Z}
%{\cal O}_{Z_{{\mathbf F}_p}})$
%is a locally free ${\cal O}_{Z_{{\mathbf F}_p}}$-module
%of rank
%${\rm codim}_XZ$
%and the sequence
%(\ref{eqXZ}) is exact except 
%possibly at 
%$F^*(N_{Z/X}\otimes_{{\cal O}_Z}
%{\cal O}_{Z_{{\mathbf F}_p}})$ by (\ref{eqFI}),
%the assertion follows.
\qed

}

\begin{cor}
Let $A$ be a regular local ring
faithfully flat
over ${\mathbf Z}_{(p)}$
and set $A_1=A/pA$.
We consider the following conditions:

{\rm (1)}
The morphism
$A_1\to F\Omega^1_A$
of $A_1$-modules
sending $1$ to
$w(p)\in F\Omega^1_A$ is a split
injection.

{\rm (2)}
$A_1$ is regular.

{\rm 1.}
We have always {\rm (1)}$\Rightarrow${\rm (2)}.

{\rm 2.}
Assume that the quotient
$A/\sqrt{pA}$ by the radical
of the principal ideal $pA$
is essentially of finite type over
a field $k_1$ with finite $p$-basis.
Then we have {\rm (2)}$\Rightarrow${\rm (1)}.
\end{cor}

\proof{
It suffices to apply Corollary
\ref{corXZx}.1
and
Corollary
\ref{corXZ} to $B=A/pA$
respectively.
\qed

}

\section{Relation with cotangent complex}

By Proposition \ref{prSA}.2,
we may sheafify the construction
of $F\Omega^1$
on a scheme $X$. 
We call $F\Omega^1_X$
the sheaf of FW-differentials on $X$.
In this section,
we study the relation of
$F\Omega^1_X$ with
cotangent complex.
Before starting,
we prepare basic properties
of sheaves of FW-differentials.

\begin{lm}\label{lmcoh}
Let $X$ be a scheme
over ${\mathbf Z}_{(p)}$.
Let
$X_{{\mathbf F}_p}$
and $F\colon X_{{\mathbf F}_p}
\to X_{{\mathbf F}_p}$
denote the closed subscheme
$X\times_{{\rm Spec}\, {\mathbf Z}}
{\rm Spec}\, {\mathbf F}_p
\subset X$
and the absolute Frobenius morphism.

{\rm 1.}
The ${\cal O}_X$-module
$F\Omega^1_X$
is a quasi-coherent ${\cal O}_
{X_{{\mathbf F}_p}}$-module.
The canonical isomorphism
{\rm (\ref{eqAB})} defines an isomorphism
\begin{equation}
F\Omega^1_X/({\cal O}_{X_{{\mathbf F}_p}}
\cdot w(p))
\to F^*\Omega^1_{X_{{\mathbf F}_p}}.
\label{eqXXp}
\end{equation}

{\rm 2.}
Assume that $X$ is noetherian
and that the reduced
part $X_{{\mathbf F}_p,{\rm red}}$ is
a scheme of finite type over
a field $k$ with finite $p$-basis.
Then, the ${\cal O}_X$-module
$F\Omega^1_X$
is a coherent ${\cal O}_{X_{{\mathbf F}_p}}$-module.
Further if $X$ is regular
of dimension $n$,
then 
$F\Omega^1_X$
is a locally free
${\cal O}_{X_{{\mathbf F}_p}}$-module
of rank $n$.
\end{lm}

\proof{
1.
If $X={\rm Spec}\, A$,
the ${\cal O}_X$-module
$F\Omega^1_X$
is defined by the $A$-module
$F\Omega^1_A$.
Hence
the ${\cal O}_X$-module
$F\Omega^1_X$
is quasi-coherent.
The ${\cal O}_X$-module
$F\Omega^1_X$
is an ${\cal O}_{X_{{\mathbf F}_p}}$-module
by Corollary \ref{corA}.1.
The isomorphism (\ref{eqXXp}) is clear from
(\ref{eqAB}).

2.
This follows from
Corollary \ref{corA}.4
and Theorem \ref{thmreg}.2.
\qed

}

\medskip
A morphism 
$f\colon X\to Y$ of schemes defines
a canonical morphism
\begin{equation}
f^*F\Omega^1_Y
\to
F\Omega^1_X
\label{eqXYsm}
\end{equation}
of ${\cal O}_X$-modules.

\medskip
We recall some of basic properties on
cotangent complexes from
\cite[Chapitres II, III]{Ill}.
For a morphism of schemes
$X\to S$, the cotangent complex
$L_{X/S}%=(L_{X/S,q},d_q)_{q\geqq 0}
$ is defined 
\cite[Chapitre II, 1.2.3]{Ill}
as a chain complex of flat
${\cal O}_X$-modules, whose cohomology
sheaves are quasi-coherent.
There is a canonical isomorphism
${\cal H}_0(L_{X/S})\to \Omega^1_{X/S}$
\cite[Chapitre II, Proposition 1.2.4.2]{Ill}.
This induces a canonical morphism
$L_{X/S}\to \Omega^1_{X/S}[0]$.

For a commutative diagram
\begin{equation}
\begin{CD}
X'@>>> S'\\
@VfVV
@VVV\\
X@>>> S,
\end{CD}
\label{eqXS}
\end{equation}
a canonical morphism
$Lf^*L_{X/S}\to L_{X'/S'}$
is defined \cite[Chapitre II, (1.2.3.2)$'$]{Ill}. 
For a morphism
$f\colon X\to Y$
of schemes over a scheme $S$,
a distinguished triangle
\begin{equation}
Lf^*L_{Y/S}\to L_{X/S}\to L_{X/Y}\to
\label{eqLXYS}
\end{equation}
is defined
\cite[Chapitre II, Proposition 2.1.2]{Ill}.

The cohomology sheaf
${\cal H}_1(L_{X/S})$
is studied as the module of imperfection
in \cite[Chapitre 0, Section 20.6]{EGA4}.
If $X\to S$ is a closed immersion 
defined by the ideal sheaf ${\cal I}_X
\subset {\cal O}_S$ and
if $N_{X/S}={\cal I}_X/{\cal I}_X^2$ denotes the conormal
sheaf,
there exists a canonical isomorphism
\begin{equation}
{\cal H}_1(L_{X/S})\to N_{X/S}
\label{eqLN}
\end{equation}
\cite[Chapitre III, Corollaire 1.2.8.1]{Ill}.
This induces a canonical morphism
$L_{X/S}\to N_{X/S}[1]$.

\begin{lm}\label{lmNOm}
{\rm 1.
(\cite[Chapitre III, Proposition 1.2.9]{Ill})}
Let $f\colon X\to Y$ be
an immersion of schemes over
a scheme $S$.
Then, the boundary morphism
$\partial\colon N_{X/Y}\to 
f^*\Omega^1_{Y/S}$
of the distinguished triangle
$Lf^*L_{Y/S}\to L_{X/S}\to L_{X/Y}
\to $ sends
$g$ to $-dg$.

{\rm 2.} {\rm (\cite[Chapitre III, 
Proposition 3.1.2 (i)$\Rightarrow$(ii)]{Ill})}
Let $X\to S$ be a smooth morphism.
Then, 
the canonical morphism
$L_{X/S}\to \Omega^1_{X/S}[0]$
is a quasi-isomorphism.

{\rm 3.} {\rm (\cite[Chapitre III, 
Proposition 3.2.4 (iii)]{Ill})}
If $X\to S$ is a regular immersion,
the canonical morphism
$L_{X/S}\to N_{X/S}[1]$
is a quasi-isomorphism.
\end{lm}

\medskip

For a scheme $E$
over ${\mathbf F}_p$,
let $F\colon E\to E=E'$
denote 
the absolute Frobenius morphism.
We canonically identify
$\Omega^1_{E/{\mathbf F}_p}
=
\Omega^1_{E/E'}$.
We study the cohomology sheaf
${\cal H}_1(L_{E/X})$
of the cotangent complex
under a certain regularity condition.

\begin{lm}\label{lmE}
Let $E$ be a scheme 
smooth over a field $k$ of
characteristic $p>0$.

{\rm 1.}
The canonical morphism
$L_{E/{\mathbf F}_p}
\to 
\Omega^1_{E/{\mathbf F}_p}[0]$
is a quasi-isomorphism
and the ${\cal O}_E$-module
$\Omega^1_{E/{\mathbf F}_p}$
is flat.

{\rm 2.}
Let $E'$ be a scheme 
smooth over a field $k'$ of
characteristic $p>0$
and 
$E'\to E$ be a morphism
of schemes.
Then, we have an exact sequence
\begin{equation}
0\to {\cal H}_1(L_{E'/E})
\to \Omega^1_{E/{\mathbf F}_p}
\otimes_{{\cal O}_E}{\cal O}_{E'}
\to 
\Omega^1_{E'/{\mathbf F}_p}
\to {\cal H}_0(L_{E'/E})\to 0
\label{eqEE'}
\end{equation}
and 
${\cal H}_q(L_{E'/E})=0$
for $q>1$.

{\rm 3. (\cite[Theorem (7.2)]{Kz})}
Let $F\colon E\to E$ denote the 
absolute Frobenius morphism.
Then, the sequence
$0\to {\cal O}_E\to F_*{\cal O}_E
\overset d
\to F_*\Omega^1_{E/{\mathbf F}_p}$
is exact.
\end{lm}

\proof{
1.
By the distinguished triangle
$L_{k/{\mathbf F}_p}
\otimes_k{\cal O}_E
\to
L_{E/{\mathbf F}_p}
\to
L_{E/k}$
and Lemma \ref{lmNOm}.2,
the assertion is reduced
to the case where
$E={\rm Spec}\, k$.
Since the formation of
cotangent complexes
commutes with limits,
we may assume $k$
is of finite type over
${\mathbf F}_p$.
Hence, we
may assume that
$k$ is the function field
of a smooth scheme
$E$ over
${\mathbf F}_p$.
Thus the assertion follows from 
Lemma \ref{lmNOm}.2.

2.
By the distinguished triangle
$L_{E/{\mathbf F}_p}
\otimes^L_{{\cal O}_E}{\cal O}_{E'}
\to
L_{E'/{\mathbf F}_p}
\to
L_{E'/E}\to $,
the assertion follows from 1
for $E$ and $E'$.

3.
We may assume that
$k$ is finitely generated
over ${\mathbf F}_p$.
Then
$k$ is isomorphic to the function field
of a scheme $S$ smooth over 
${\mathbf F}_p$.
We may assume that
$E$ is the generic fiber
of a smooth scheme
$E_S$ over $S$.
Thus, it is reduced to
the case where $k
={\mathbf F}_p$
is perfect.
Then, 
the canonical morphism
$\Omega^1_{E/{\mathbf F}_p}
\to \Omega^1_{E/k}$
is an isomorphism
and the assertion follows from 
the Cartier isomorphism
\cite[Theorem (7.2)]{Kz}.
\qed

}

\begin{lm}\label{lmEX}
Let $X$ be a scheme.
Let $p$ be a prime number and
$E$ be a scheme over
${\mathbf F}_p$.
Let $f\colon E\to X$ be a morphism 
of schemes.

{\rm 1.}
We consider the following conditions:

{\rm (1)}
The morphism
$f\colon E\to X$ factors through
the absolute Frobenius morphism
$F\colon E\to E$.

{\rm (2)}
The canonical surjection
\begin{equation}
%\Omega^1_{E/k}
%=
\Omega^1_{E/{\mathbf F}_p}=
\Omega^1_{E/{\mathbf Z}}
\to
\Omega^1_{E/X}
\label{lmOmE}
\end{equation}
is an isomorphism.

%\noindent
We have
{\rm (1)}$\Rightarrow${\rm (2)}.
If $E$ is a smooth scheme
over a field $k$,
we have
{\rm (2)}$\Rightarrow${\rm (1)}.

{\rm 2.}
Assume that $X$ is a regular noetherian
scheme,
that $E$ is smooth over
a field
and that $f$ is of finite type
and satisfies
the equivalent conditions in {\rm 1.}
Then the ${\cal O}_E$-module
${\cal H}_1(L_{E/X})$ is locally free
of finite rank.
\end{lm}

\proof{
1.
(1)$\Rightarrow$(2):
Suppose $f\colon E\to X$ factors through
$F\colon E\to E=E'$.
Then since the surjection
$\Omega^1_{E/{\mathbf F}_p}
\to 
\Omega^1_{E/E'}$
is an isomorphism,
the surjections
$\Omega^1_{E/{\mathbf Z}}\to
\Omega^1_{E/X}
\to 
\Omega^1_{E/E'}$
are isomorphisms.

(2)$\Rightarrow$(1):
The condition (2)
means that
the composition of
$f^{-1}{\cal O}_X\to {\cal O}_E$
and $d\colon  {\cal O}_E\to
\Omega^1_{E/{\mathbf F}_p}$ is the $0$-morphism.
Since $F\colon E\to E$
is a homeomorphism
on the underlying topological
spaces,
the continuous mapping
$f\colon E\to X$
is the composition of
$F\colon E\to E$
with a unique continuous mapping
$g\colon E\to X$.
Thus, the condition
(2) is equivalent
to the condition that the composition
$g^{-1}{\cal O}_X\to
F_*{\cal O}_E
\to F_*\Omega^1_{E/{\mathbf F}_p}$ 
is the $0$-morphism.

By Lemma \ref{lmE}.3,
the sequence
$0\to {\cal O}_{E}
\to F_*{\cal O}_E
\overset{d}\to F_*\Omega^1_{E/
{\mathbf F}_p}$
is exact.
%means isomorphism
%$\Omega^1_{E/s}\to
%\Omega^1_{E/X}$
Thus, the condition
(2) is further equivalent
to the condition that 
the morphism
$g^{-1}{\cal O}_X\to
F_*{\cal O}_E$ factors through
$g^{-1}{\cal O}_X\to
{\cal O}_E$.
Since $F\colon E\to E$ is
affine, this defines
a morphism $g\colon E\to X$
of schemes and
the condition (2) 
is equivalent to (1).

2.
Since the assertion is local on $E$,
we may assume that 
$E$ and $X$ are affine
and there exists
a closed immersion
$E\to P={\mathbf A}^n_X$
for some $n$.
Since $E$ and $X$
hence $P$ are regular,
the closed immersion
$E\to P$ is a regular
immersion.
Then, the distinguished triangle
$L_{P/X}\otimes_{{\cal O}_P}
{\cal O}_E\to
L_{E/X}\to L_{E/P}\to$ (\ref{eqLXYS})
defines an exact sequence
$0\to
{\cal H}_1(L_{E/X})\to N_{E/P}\to
\Omega^1_{P/X}\otimes_{{\cal O}_P}
{\cal O}_E
\to 
\Omega^1_{E/X}\to 0$
by Lemma \ref{lmNOm}.2
for $P\to X$
and Lemma \ref{lmNOm}.3
for $E\to P$.
The ${\cal O}_E$-modules
in the exact sequence
other than $
{\cal H}_1(L_{E/X})$ are  locally free
of finite rank 
by the isomorphism (\ref{lmOmE}).
Hence
${\cal H}_1(L_{E/X})$ is also locally free
of finite rank.
\qed

}
\medskip

We give a constuction
yielding an FW-derivation.

\begin{lm}\label{lmdu}
Let $X$ be a scheme
and set ${\mathbf A}^1_X
=X\times_{
{\rm Spec}\, {\mathbf Z}}
{\rm Spec}\, {\mathbf Z}[T]$.

{\rm 1.}
Let $E$ be a scheme 
over ${\mathbf F}_p$
and let $E\to {\mathbf A}^1_X$
be a morphism of schemes.
Then, the
distinguished triangle
$
L_{{\mathbf A}^1_X/X}
\otimes^L_{{\cal O}_{{\mathbf A}^1_X}}
{\cal O}_E
\to
L_{E/X}
\to
L_{E/{\mathbf A}^1_X}
\to $
defines an exact sequence
\begin{equation}\begin{CD}
0@>>>
{\cal H}_1(L_{E/X})
@>>>
{\cal H}_1(L_{E/{\mathbf A}^1_X})
@>>>
\Omega^1_{{\mathbf A}^1_X/X}
\otimes_{{\cal O}_{ {\mathbf A}^1_X}}
{\cal O}_E.
\end{CD}
\label{eqLW}
\end{equation}

{\rm 2.}
Let $u\in \Gamma(X,{\cal O}_X)$.
Define a closed subscheme
$W\subset {\mathbf A}^1_X$
by the ideal $(u-T^p,p)$
and identify
${\cal H}_1(L_{W/{\mathbf A}^1_X})$
with the conormal sheaf
$N_{W/{\mathbf A}^1_X}$
by the canonical isomorphism
{\rm (\ref{eqLN})}.
Then, the section
$u-T^p$ of the conormal sheaf
$N_{W/{\mathbf A}^1_X}$ lies in the image
of the injection
\begin{equation}
\Gamma(W,{\cal H}_1(L_{W/X}))
\to
\Gamma(W,{\cal H}_1(L_{W/ {\mathbf A}^1_X}))
=
\Gamma(W,N_{W/ {\mathbf A}^1_X})
\label{eqWA}
\end{equation}
defined by {\rm (\ref{eqLW})}
for $E=W$.
In other words,
there exists a unique section
\begin{equation}
 \omega
\in \Gamma(W,{\cal H}_1(L_{W/X}))
\label{eqtilom}
\end{equation}
such that the image
in $\Gamma(W,N_{W/ {\mathbf A}^1_X})$
equals
$u-T^p$.
\end{lm}

\proof{
1.
Since the ${\cal O}_{{\mathbf A}^1_X}$-module $\Omega^1_{{\mathbf A}^1_X/X}$
is flat,
the assertion follows from
the canonical isomorphism
$L_{{\mathbf A}^1_X/X}
\to
\Omega^1_{{\mathbf A}^1_X/X}[0]
$ in Lemma \ref{lmNOm}.2.

2.
By 1 applied to $E=W$,
to show that $u-T^p$ lies in the image
of (\ref{eqWA}),
it suffices to show that this vanishes in
$\Gamma(W,
\Omega^1_{{\mathbf A}^1_X/X}
\otimes_{{\cal O}_{{\mathbf A}^1_X}}
{\cal O}_W)$.
By Lemma \ref{lmNOm}.1,
the last arrow in (\ref{eqLW})
for $E=W$
is $-d\colon
N_{W/ {\mathbf A}^1_X}\to
\Omega^1_{{\mathbf A}^1_X/X}
\otimes_{{\cal O}_{{\mathbf A}^1_X}}
{\cal O}_W$.
Since
$d(u-T^p)=-pT^{p-1}dT=0$
on $W$,
the assertion follows.
\qed

}

\begin{df}[{\rm cf.\ \cite[
Definition 1.1.6]{gr}
or \cite[Definition 1.1.6 in v1]{lex}}] \label{dfw}
Let $X$ be a scheme
and $u\in \Gamma(X,{\cal O}_X)$
be a section.
Let $E$ be a scheme over
${\mathbf F}_p$
and 
let $f\colon E\to X$ be a morphism
of schemes.
Let $v\in \Gamma(E,{\cal O}_E)$
be a section such that
$u|_E=f^*u\in  \Gamma(E,{\cal O}_E)$
is the $p$-th power of $v$.
Let $W\subset {\mathbf A}^1_X$
be the closed subscheme
as in Lemma {\rm \ref{lmdu}}
and 
define a morphism
$E\to W$ over $X$
by sending $T$ to 
$v\in \Gamma(E,{\cal O}_E)$.
We define a section
\begin{equation}
w(u,v)\in \Gamma(E,{\cal H}_1(L_{E/X}))
\label{eqomega}
\end{equation}
to be the image of
$\omega$ in {\rm (\ref{eqtilom})}
by the morphism
$\Gamma(W,{\cal H}_1(L_{W/X}))
\to
\Gamma(E,{\cal H}_1(L_{E/X}))$
defined by $E\to W$.
\end{df}

\begin{pr}[{\rm cf.\ \cite[Lemma 1.1.4]{gr}
or \cite[Proposition 1.1.5 in v1]{lex}
}] \label{prdu}
Let $X$ be a scheme
and $u\in \Gamma(X,{\cal O}_X)$.
Let $f\colon E\to X$ be a morphism
of schemes
and assume that $E$
is a scheme over ${\mathbf F}_p$.
Let $v\in \Gamma(E,{\cal O}_E)$
be a section satisfying
$u|_E=f^*u\in  \Gamma(E,{\cal O}_E)$
is the $p$-th power of $v$.

{\rm 1.}
Assume $u|_E=0$
and let $E\to Z\subset X$
be the morphism
to the closed subscheme 
defined by $u$.
Then 
$w(u,0)\in \Gamma(E,{\cal H}_1(L_{E/X}))$
is the image of
$u\in \Gamma(Z,N_{Z/X})$
by the morphism
$\Gamma(Z,N_{Z/X})
\to \Gamma(E,{\cal H}_1(L_{E/X}))$
defined
by $L_{Z/X}\otimes^L
_{{\cal O}_Z}{\cal O}_E
\to L_{E/X}$.

{\rm 2.}
Let
$u'\in \Gamma(X,{\cal O}_X)$
and $v'\in \Gamma(E,{\cal O}_E)$
be another pair of sections
satisfying
$u'|_E=v'^p$.
Then,
we have
\begin{align}
w(u+u',v+v')
&\, =
w(u,v)+
w(u',v')-
P(v,v')
\cdot w(p,0),
\label{eqaddL}
\\
w(uu',vv')
&\, =
u'\cdot w(u,v)+
u\cdot w(u',v').
\label{eqLbL}
\end{align}

{\rm 3.}
Let $X\to S$ be a morphism
of schemes.
Then, the minus of the boundary mapping
$-\partial\colon {\cal H}_1(L_{E/X})
\to \Omega^1_{X/S}
\otimes_{{\cal O}_X}{\cal O}_E$
of the distinguished triangle $L_{X/S}
\otimes_{{\cal O}_X}^L{\cal O}_E
\to
L_{E/S}\to L_{E/X}\to$
sends
$w(u,v)\in \Gamma(E,{\cal H}_1(L_{E/X}))$
to $
du\in \Gamma(E,\Omega^1_{X/S}
\otimes_{{\cal O}_X}{\cal O}_E)$.

\end{pr}

%\begin{pr}[{cf.\ \cite[]{gr}}]\label{prdu2}
%
%\end{pr}

\proof{
1.
Since the morphism
$E\to W\subset 
{\mathbf A}^1_X$ factors through
the $0$-section $Z\subset 
{\mathbf A}^1_X$,
the assertion follows from
$T^p=0$
in $\Gamma(Z,N_{Z/{\mathbf A}^1_X})$.

2.
By 1, 
$w(p,0)\in\Gamma(E,{\cal H}_1(L_{E/X}))$
is the image
of $p\in N_{{\mathbf F}_p/
{\mathbf Z}}$.
Let $W'$ be the
closed subscheme of
${\mathbf A}^2_X$
defined by the ideal $(T^p-u,T'^p-u',p)$
and define $E\to W'$
by $T\mapsto v$, $T'\mapsto v'$.
Then,
(\ref{eqaddL})
follows from the binomial expansion
$$(u+u')-(T+T')^p
=
(u-T^p)
+(u'-T^{\prime p})
-
P(T,T')\cdot p$$
Similarly,
(\ref{eqLbL})
follows from
$$(uu')-(TT')^p
=
u'(u-T^p)
+u(u'-T^{\prime p})
-
(u-T^p)
(u'-T^{\prime p}).$$

3.
The morphisms
$E\to W\to {\mathbf A}^1_X\to X\to S$
define a commutative diagram
$$\begin{CD}
{\cal H}_1(L_{E/X})
@>>>
{\cal H}_1(L_{E/{\mathbf A}^1_X})
@<<<
N_{W/{\mathbf A}^1_X}
\otimes_{{\cal O}_W}{\cal O}_E
\\
@V{-\partial}VV@V{-\partial}VV@VVdV\\
\Omega^1_{X/S}
\otimes_{{\cal O}_X}{\cal O}_E
@>>>
\Omega^1_{{\mathbf A}^1_X/S}
\otimes_{{\cal O}_{{\mathbf A}^1_X}}
{\cal O}_E
@=
\Omega^1_{{\mathbf A}^1_X/S}
\otimes_{{\cal O}_{{\mathbf A}^1_X}}
{\cal O}_E
\end{CD}$$
by Lemma \ref{lmNOm}.1.
Since $d(u-T^p)=du$
in $\Gamma(E,
\Omega^1_{{\mathbf A}^1_X/S}
\otimes_{{\cal O}_{{\mathbf A}^1_X}}
{\cal O}_E)$
and since the lower left
horizontal arrow is an injection,
the assertion follows.
\qed

}

\begin{cor}\label{cordu}
Let $X$ be a scheme
and let $E$ be a scheme over
${\mathbf F}_p$.
Let $g\colon E\to X$ be
a morphism of schemes
and let $L_{E/X}$
denote the cotangent
complex for the composition $f=g\circ F
\colon E\to X$ with
the absolute Frobenius $F\colon
E\to E$.
Then, the mapping
\begin{equation}
w\colon
\Gamma(X,{\cal O}_X)
\to \Gamma(E,{\cal H}_1(L_{E/X}))
\label{eqdfdu}
\end{equation}
sending 
$u\in \Gamma(X,{\cal O}_X)$
to 
$w(u,v)$
for $v=g^*u\in \Gamma(E,{\cal O}_E)$
is an FW-derivation.
\end{cor}

\proof{
The assertion follows
from Proposition \ref{prdu}.2.
\qed

}
\medskip

The construction of 
the FW-derivation
$w$ (\ref{eqdfdu}) is functorial
in $X$ and $E$.

\begin{df}\label{dfdu}
Let $X$ be a scheme
and let $E$ be a scheme over
${\mathbf F}_p$.
Let $g\colon E\to X$ be
a morphism of schemes
and let $L_{E/X}$
denote the cotangent
complex for the composition $f=g\circ F
\colon E\to X$ with
the absolute Frobenius $F\colon
E\to E$.
By sheafifying the morphism
{\rm  (\ref{eqdfdu})},
we define an FW-derivation
$w\colon
g^{-1}{\cal O}_X\to
{\cal H}_1(L_{E/X})$
and the morphism
\begin{equation}
g^*F\Omega^1_X
\to {\cal H}_1(L_{E/X})
\label{eqwd}
\end{equation}
defined by the universality
of $F\Omega^1_X$.
\end{df}
\medskip

We study condition for
the morphism
(\ref{eqwd}) to be an isomorphism.

\begin{lm}\label{lmdu2}
Let $g\colon E\to Z$ be
a morphism of schemes
over
${\mathbf F}_p$
and and let $L_{E/Z}$
denote the cotangent
complex for the composition $f=g\circ F
\colon E\to Z$ with
the absolute Frobenius $F\colon
E\to E$.

{\rm 1.}
The morphism
$g^*F\Omega^1_Z
\to {\cal H}_1(L_{E/Z})$
{\rm (\ref{eqwd})}
is a split injection.

{\rm 2.}
The split injection
{\rm (\ref{eqwd})}
is an isomorphism if 
${\cal H}_1(L_{E/{\mathbf F}_p})=0$.
The condition 
${\cal H}_1(L_{E/{\mathbf F}_p})=0$
is satisfied if
$E$ is smooth over a field.
\end{lm}

\proof{
1.
The composition
$$\begin{CD}
g^*F\Omega^1_Z
@>{\rm (\ref{eqwd})}>> 
{\cal H}_1(L_{E/Z})
@>{-\partial}>>
f^*\Omega^1_{Z/{\mathbf F}_p}
\end{CD}$$
is the isomorphism induced
by (\ref{eqB})
by Proposition \ref{prdu}.3.
Hence 
$g^*F\Omega^1_Z
\to {\cal H}_1(L_{Z/X})$
{\rm (\ref{eqwd})}
is a split injection.

2.
The distinguished triangle
$Lf^*L_{Z/{\mathbf F}_p}
\to
L_{E/{\mathbf F}_p}\to L_{E/Z}
\to$
defines an exact sequence
$
{\cal H}_1(L_{E/{\mathbf F}_p})
\to
{\cal H}_1(L_{E/Z})
\to
f^*\Omega^1_{Z/{\mathbf F}_p}$.
Hence the vanishing
${\cal H}_1(L_{E/{\mathbf F}_p})=0$
implies the isomorphism.

If $E$ is smooth over a field,
we have
${\cal H}_1(L_{E/{\mathbf F}_p})=0$
by Lemma \ref{lmE}.1.
\qed

}

\begin{pr}\label{prdu2}
Let $X$ be a scheme
and let $E$ be a scheme over
${\mathbf F}_p$.
Let $g\colon E\to X$ be
a morphism of schemes and
$Z\subset X$ be a closed
subscheme such that
$g\colon E\to X$ factors through
$g_Z\colon E\to Z$ and that
$Z$ is a scheme over ${\mathbf F}_p$.
Let $L_{E/X}$ and $L_{E/Z}$
denote the cotangent
complexes for the compositions
$f=g\circ F
\colon E\to X$ 
and $f_Z=g_Z\circ F\colon E\to Z$
with the absolute Frobenius $F\colon
E\to E$.

{\rm 1.}
The canonical morphism
$g^*F\Omega^1_X
\to {\cal H}_1(L_{E/X})$
{\rm (\ref{eqwd})}
is a surjection
if ${\cal H}_1(L_{E/{\mathbf F}_p})
=0$.
The condition 
${\cal H}_1(L_{E/{\mathbf F}_p})
=0$
is satisfied if $E$ is smooth over 
a field.

{\rm 2.}
The canonical morphism
$g^*F\Omega^1_X
\to {\cal H}_1(L_{E/X})$
{\rm (\ref{eqwd})}
and the morphism
$f_Z^*N_{Z/X}
\to 
g^*F\Omega^1_X$
defined by
{\rm (\ref{eqFI})}
are injections
if ${\cal H}_2(L_{E/Z})=0$.

The condition ${\cal H}_2(L_{E/Z})=0$
is satisfied if $E$ and $Z$
are smooth over fields.
\end{pr}

\proof{
We consider the commutative diagram
\begin{equation}
\begin{CD}
@.f_Z^*N_{Z/X}
@>>> 
g^*F\Omega^1_X
@>>> 
g_Z^*F\Omega^1_Z
@>>>0
\\
@.@|@V{\rm (\ref{eqwd})}VV
@V{\rm (\ref{eqwd})}VV
@.\\
{\cal H}_2(L_{E/Z})
@>>>
f_Z^*N_{Z/X}
@>>>
{\cal H}_1(L_{E/X})
@>>>
{\cal H}_1(L_{E/Z})
@>>>0
\end{CD}
\label{eqOmL}
\end{equation}
of exact sequences.
The lower line is defined
by the distinguished triangle
$Lf_Z^*L_{Z/X}
\to
L_{E/X}\to L_{E/Z}
\to$
and 
the upper line is the pull-back of
the exact sequence
defined by (\ref{eqFI}).

1.
If ${\cal H}_1(L_{E/{\mathbf F}_p})=0$,
the right vertical arrow is
an isomorphism by Lemma \ref{lmdu2}.
Hence the middle vertical arrow
is a surjection.
If $E$ is smooth over a field,
we have
${\cal H}_1(L_{E/{\mathbf F}_p})=0$
by Lemma \ref{lmE}.1.

2.
If ${\cal H}_2(L_{E/Z})=0$,
since the right vertical arrow is
an injection by Lemma \ref{lmdu2},
the middle vertical arrow is
an injection.
Further the morphism
$f_Z^*N_{Z/X}
\to
g^*F\Omega^1_X$
is an injection
by the commutativity of
the left square.

If $E$ and $Z$
are smooth
over fields,
we have
${\cal H}_2(L_{E/Z})=0$
by Lemma \ref{lmE}.2.
\qed

}

\begin{cor}\label{cordu2}
Let $A$ be a local ring with residue field
$k$ of characteristic $p>0$.
Then, the canonical morphism
$F\Omega^1_A\otimes_Ak
\to 
{\cal H}_1(L_{k/A})$ 
{\rm(\ref{eqwd})}
is an isomorphism.
\end{cor}

\proof{
It suffices to apply Proposition \ref{prdu2}
to $g\colon Z={\rm Spec}\, k
\to X={\rm Spec}\, A$.
\qed

}

\end{document}